\documentclass[11pt,english]{article}

\usepackage{a4wide}
\usepackage{amssymb,amsthm}
\usepackage{amsmath}
\usepackage{amsfonts}
\usepackage{latexsym}
\usepackage{amsxtra}
\usepackage{graphicx}
\usepackage{wrapfig}
\usepackage{mathrsfs}
\usepackage{epsfig}
\usepackage{times}
\usepackage{makeidx}
\usepackage[titles]{tocloft}
\usepackage{cite}

\DeclareMathAlphabet{\mathpzc}{OT1}{pzc}{m}{it}

\newcommand{\heikodetail}[1]{}
\renewcommand{\proof}{\noindent{\sc Proof:}\hskip 1.1em}
\renewcommand{\qed}{\hfill\mbox{$\Box$}\\}
\def\yy{\mbox{$\spadesuit$}}
\newcommand{\invisible}[1]{\par ($\spadesuit$ \emph{hidden comments in \TeX{}
file\/})} 
\newcommand{\reallyinvisible}[1]{}      

\newfont{\thickmath}{msbm10 scaled \magstephalf}%
\newfont{\smallthickmath}{msbm7 scaled \magstephalf}%
\newfont{\footnotethickmath}{msbm8}%
\newfont{\footnotesmallthickmath}{msbm6}%

\newcommand{\cA}{{\mathcal{A}}}

\newcommand{\cC}{{\mathcal{C}}}

\newcommand{\cR}{{\mathcal{R}}}
\newcommand{\cT}{{\mathcal{T}}}
\newcommand{\area}{{\textnormal{area}}} 
\newcommand{\sym}{{\textnormal{sym}}} 
\newcommand{\vol}{{\textnormal{vol}}} 

\newcommand{\hR}{{\widehat{\mathcal{R}}}}
\newcommand{\hT}{{\widehat{\mathcal{T}}}}
\newtheorem{lemma}{\bf Lemma}[section]
\newtheorem{theorem}[lemma]{\bf Theorem}

\newtheorem{corollary}[lemma]{\bf Corollary}

\newtheorem{definition}[lemma]{\bf Definition}

\newcommand{\Fo}{\,\,\,\text{for }\,\,}
\newcommand{\Foa}{\,\,\,\text{for all }\,\,}

\newcommand{\AND}{\,\,\,\text{and }\,\,}

\newcommand\Reals{{\mathbb R}}
\newcommand\R{{\mathbb R}}

\newcommand\N{{\mathbb N}}
\renewcommand\S{{\mathbb S}}

\newcommand{\bbbr}{\Reals}

\newcommand\dist{\mathop{\rm dist}\nolimits}
\newcommand\Int{\mathop{\rm int}\nolimits}

\newcommand\Id{{{\rm Id}}}

\newcommand{\xx}{\mbox{$\clubsuit$}}

\renewcommand{\R}{\bbbr}

\newcommand{\eps}{\varepsilon}
\renewcommand{\H}{\mathscr{H}}

\begin{document}

\renewcommand{\thefigure}{\arabic{figure}}

\title{\large\bf On minimal immersions in Finsler space.}

\author{\normalsize Patrick Overath, Heiko von der Mosel}


\maketitle


\frenchspacing

\begin{abstract}
We explore a connection between the Finslerian area functional and
well-investigated Cartan functionals to prove new Bernstein theorems, 
uniqueness and removability results for Finsler-minimal graphs, as well as enclosure theorems
and isoperimetric inequalities for minimal immersions  in Finsler
spaces. In addition, we establish the existence of smooth Finsler-minimal
immersions spanning given extreme or graphlike boundary contours.

\vspace{2mm}


\end{abstract}

      


\renewcommand\theequation{{\thesection{}.\arabic{equation}}}
\def\setnumbers{\setcounter{equation}{0}}


\section{Introduction}\label{sec:1}
\subsection{Minimal immersions in Finsler geometry}
\label{sec:1.1}
Minimal surface theory in Finsler spaces seems to be a largely underdeveloped
terrain. Recently, we established a new connection between Finsler-minimal immersions and anisotropic variational integrals, so-called {\it Cartan functionals}, to treat
the Plateau problem in Finsler $3$-space; see  \cite{overath-vdM-2012a}.
In the present note we explore this connection  to substantially extend the few 
known results about minimal graphs, available only in very specific Finsler-Minkowski spaces so far,
and we prove new global results for
Finsler-minimal immersions such as enclosure theorems
and isoperimetric inequalities. In addition, we establish the existence of smooth
Finsler-minimal immersions spanning given extreme or graphlike boundary contours.
Our general assumption on the Finsler structure
turns out to be natural and sharp, since it translates to known sharp threshold
values for the anisotropies of the few specific Finsler-Minkowski spaces investigated
so far.

\medskip

In order to
 briefly recall the precise 
notion of Finsler-minimal immersions (in the sense of Busemann and Hausdorff) let
$\mathscr{N}=\mathscr{N}^n$ be an $n$-dimensional smooth
manifold with tangent bundle $T\mathscr{N}:=\bigcup_{x\in\mathscr{N}}T_x\mathscr{N}$ and its zero-section $o:=\{(x,0)\in T\mathscr{N}\}$.
A non-negative function $F\in C^\infty(T\mathscr{N}\setminus o)$ is called
a {\it Finsler metric} on $\mathscr{N}$ (so that $(\mathscr{N},F)$ becomes
a {\it Finsler manifold}) if $F$ satisfies the conditions
\begin{enumerate}
\item[\rm (F1)] $F(x,ty)=tF(x,y)$ for all $t>0$ and all $(x,y)\in T\mathscr{N}$ (homogeneity);
\item[\rm (F2)] $g_{ij}(x,y):=\big(F^2/2)_{y^iy^j}(x,y)$ form the coefficients
of a positive definite matrix, the {\it fundamental tensor,}
for all $(x,y)\in
T\mathscr{N}\setminus o $, where  for  given local coordinates $x^1,\ldots,x^n$ about
$x\in \mathscr{N}$, the $y^i,$ $i=1,\ldots,n$, denote the corresponding
bundle coordinates via $y=y^i\frac{\partial}{\partial  x^i}|_x\in 
T_x\mathscr{N}.$ Here we sum over repeated Latin indices from $1$ to $n$
according to the Einstein summation convention, and
$F(x,y)$ is written as $F(x^1,\ldots,x^n,y^1,\ldots,y^n)$.
\end{enumerate}
If $F(x,y)=F(x,-y)$ for all $(x,y)\in T\mathscr{N}$ then $F$ is
called a {\it reversible} Finsler metric, 
\heikodetail{

\bigskip

\tt ICH HATTE GESCHRIEBEN: and if there are coordinates such that
$F$ depends only on $y$, then $F$ is called a {\it Minkowski metric}.

Herr Overaths Antwort: \yy\yy  Wir brauchen einen Minkowski space, d.h. einen 
reellen Vektorraum auf welchem die Finsler metric invariant unter 
translation ist. Kurz einen $\R^n$ auf welchem $F=F(y)$ in 
Standardkoordinaten (vgl. \cite{bao-chern-shen_2000}, 
\cite[Definition 1.2.1, p. 6]{shen-lecture-2001}). Meiner Meinung nach 
stellt das einen Unterschied dazu dar, spezielle andere Koordinaten 
vorauszusetzen in denen $F=F(y)$ gilt. Ich schlage vor:

\bigskip

} 
and if $F$ on $\mathscr{N}=\R^n$ depends only on $y$ in standard coordinates, then $F$ is called a {\it Minkowski metric}.
\heikodetail{

\bigskip

Wichtige Bemerkung von Herrn Overath:
{\tt \yy In the theorems, we use Minkowski metrics defined on a vector space, so rather choose $F$ to be called a Minkowski metric, if it depends only on $y$ w.r.t. coordinates compatible with the vector space structure. Remember, that we use theorems on parametric Lagrangians proven in a vector space ($\R^{n}$) setting. }

\bigskip

}
Any $C^2$-immersion $X:\mathscr{M}^m\hookrightarrow \mathscr{N}^n$ from a
smooth 
$m$-dimensional manifold $\mathscr{M}=\mathscr{M}^m$
into
$\mathscr{N}$ induces a {\it pulled-back Finsler metric} $X^*F$ on $\mathscr{M}$
via
$$
(X^*F)(u,v):=F(X(u),dX|_u(v))\quad\textnormal{for $(u,v)\in T\mathscr{M}.$}
$$
Following Busemann \cite{busemann-1947} and Shen \cite{shen98} we define
the {\it Busemann-Hausdorff volume form} as the volume ratio of the
Euclidean and the Finslerian unit ball, i.e.,
$$
d\textnormal{vol}_{X^*F}(u):=\sigma_{X^*F}(u)du^1\wedge \ldots\wedge du^m
\quad\textnormal{on $\mathscr{M}$,}
$$
where
\begin{equation}\label{volumefactor}
\sigma_{X^*F}(u):=\frac{\mathscr{H}^m(B_1^m(0))}{\mathscr{H}^m(
\{v=(v^1,\ldots, v^m)\in\R^m:X^*F(u,v^\delta\frac{\partial}{\partial u^\delta|_u})\le 1\}},
\end{equation}
with a summation over  Greek indices from $1$ to $m$ in the denominator.
Here $\mathscr{H}^m$ denotes the $m$-dimensional Hausdorff-measure.
The {\it Busemann-Hausdorff area} or in short {\it Finsler area}\footnote{Notice that the alternative Holmes-Thompson volume form (see \cite{alvarez-berck-wrong-hausdorff-2006}) leads to a different notion of Finslerian minimal surfaces that we do not
address here.} of the immersion $X:\mathscr{M}\to\mathscr{N}$ is then given
by
\begin{equation}\label{area}
\area^F_\Omega(X):=\int_{u\in\Omega}\,d\vol_{X^*F}(u)
\end{equation}
for a measurable subset $\Omega\subset\mathscr{M}.$
Shen \cite[Theorem 1.2]{shen98}  derived
the first variation of this functional which leads to the definition
of Finsler-mean curvature. Critical immersions for $\area^F_\Omega$
are therefore {\it Finsler-minimal immersions}, or simply minimal 
hypersurfaces 
in $(\mathscr{N},F)$ if the co-dimension $n-m$ equals $1$.

By means of his variational formulas Shen excluded the existence of
closed oriented Finsler-minimal submanifolds in Minkowski space, i.e., in $\R^n$ equipped
with a Minkowski metric $F=F(y)$; see \cite[Theorem 1.3]{shen98}.
 Souza and Tenenblatt proved in 
\cite{souza-tenenblat-2003} that there is  up to homotheties exactly one
complete embedded Finsler-minimal surface of revolution in the
{\it Minkowski-Randers space} $(\R^3,F)$ where $F(y):=|y|+b_iy^i$ for
some constant vector $b$ with length $|b|<1.$ In the same setting
Souza, Spruck, and Tenenblatt \cite{sst} computed the fairly complicated
pde for Finsler-minimal graphs, and proved under  the more restrictive bound $|b|<1/\sqrt{3}$ that any  Finsler-minimal graph defined on the entire plane $\R^2$
is an affine plane. This bound on $|b|$ is actually sharp for
this Bernstein theorem: beyond the threshold value, i.e., for $|b|\in (1/\sqrt{3},1)$,
where $(\R^3,F)$ is still a Finsler space (see e.g. \cite[p. 4]{chern-shen-2005}),
 the pde
ceases to be elliptic, and there exists 
a Finsler-minimal cone with a point singularity in that case; see
\cite[Proposition 13]{souza-tenenblat-2003} and \cite[p. 300]{sst}.  
Souza et al. 
also proved a touching principle and the removability of isolated singularities
of Finsler-minimal graphs in this specific three-dimensional Minkowski-Randers space.
Their Bernstein theorem was later generalized by Cui and Shen \cite{cui-shen-2009}
to entire
$m$-dimensional
Finsler-minimal graphs in the more general {\it $(\alpha,\beta)$-Minkowski spaces} 
$(\R^{m+1},F) $. Here,  $F(y)$ equals $\alpha(y)\phi\big[\beta(y)/\alpha(y)\big]$
with $\alpha(y):=|y|$ and  the linear perturbation term
$\beta(y):=b_iy^i$, where $\phi$ is a positive smooth scalar function
satisfying a particular differential equation to guarantee
that $F$ is at least a Finsler metric; see e.g. \cite[Lemma 1.1.2]{chern-shen-2005}.
Cui and Shen require for their Bernstein results 
fairly complicated additional and more restrictive
conditions on $\phi$ (see condition 
(1) in \cite[Theorem 1.1]{cui-shen-2009} or condition (4) of
\cite[Theorem 1.2]{cui-shen-2009}) that could be verified only for
a few 
specific choices of $(\alpha,\beta)$-metrics, and only in dimension $m=2$: 
for the Minkowski-Randers case with
$\phi(s)=1+s$ if $|b|<1/\sqrt{3}$ (reproducing \cite[Theorem 6]{sst}),
 for the {\it two-order metric} with $\phi(s)=(1+s)^2$ 
under the condition $|b|<1/\sqrt{10}$, or for the {\it Matsumoto metric}
where $\phi(s)=(1-s)^{-1}$ if $|b|<1/2.$ Also these threshold values
for $|b|$ are sharp; one finds Finsler-minimal cones for $|b|$ beyond the respective bounds for the Minkowski-Randers and for the two-order metric. And the
Matsumoto metric simply ceases to be a Finsler metric if $|b|>1/2.$; see \cite[Theorem 5.3 \& Section 6]{cui-shen-2009}. 
Cui and Shen \cite{cui-shen-2011}
 derived a representation formula for rotationally
symmetric Finsler-minimal surfaces in Minkowski-$(\alpha,\beta)$-spaces 
\heikodetail{

\bigskip

Herr Overath:
In \cite{he-yang} wird wohl 
\"Ahnliches aber f\"ur nur das Holmes-Thompson-Volumen 
betrachtet so weit ich das beurteilen kann!

\bigskip

}
and they  presented a unique
explicit forward-complete rotationally symmetric Finsler-minimal surface in Minkowski-Randers-$3$-space. 

Our results presented in the next subsection generalize and extend many
of these results to
general Finsler-Minkowski spaces 
under a natural assumption on a suitable symmetrization of
the underlying Finsler metric.
It turns out that this assumption reproduces the sharp bounds on the anisotropy
$|b|$ in the specific Finsler-Minkowski-spaces described above.

For any function $F\in C^0(T\mathscr{N}\setminus o)$ that is positively $1$-homogeneous
in the $y$-variable, we define the
{$m$-harmonic symmetrization} $F_\sym$ as
\begin{equation}\label{m-harmonic}
F_{\textnormal{sym}}(x,y):= \left[\frac{2}{\frac{1}{F^m(x,y)}+\frac{1}{F^m(x,-y)}}\right]^{\frac{1}{m}}\quad\Fo (x,y)\in T\mathscr{N}\setminus o,
\end{equation}
which by definition is even and positively $1$-homogeneous in the $y$-variable,
and thus
continuously extendible by zero to the whole tangent bundle $T\mathscr{N}$.
But even if $F$ is a Finsler metric, the symmetrized form $F_\sym$ might not be,
which motivates our {\bf General Assumption}:
\begin{enumerate}
\item[\rm\bf (GA)]
{\it Let $F(x,y)$ be a Finsler metric on $\mathscr{N}=\R^{m+1}$ (with respect to its
standard coordinates) such that its
$m$-harmonic symmetrization
$F_\sym(x,y)$ is also a Finsler metric on $\R^{m+1}$.}
\end{enumerate}
Notice that a reversible Finsler metric $F$ automatically coincides with
its $m$-harmonic symmetrization $F_\sym$ so that our general assumption
(GA) is superfluous in reversible Finsler spaces.
A sufficient criterion guaranteeing that (GA) holds for non-reversible Finsler
metrics was derived in \cite[Theorem 1.5]{overath-vdM-2012a}, allowing for
a non-trivial $x$-dependence such as $F(x,y):=F_\textnormal{rev}(x,y)+b_iy^i$ where 
$F_\textnormal{rev}$ is
a reversible Finsler metric. Even in the Minkowski setting this creates
examples that were not treated before, for example, choosing the {\it perturbed
quartic metric} (see \cite[p. 15]{bao-chern-shen_2000})
$$
F_\textnormal{rev}(y):=\sqrt{\sqrt{\sum_{i=1}^{m+1} (y^i)^4} + 
\epsilon \sum_{i=1}^{m+1} (y^i)^2}\quad\Fo \epsilon >0
$$
as the reversible part of $F(y)=F_\textnormal{rev}(y)+b_iy^i$.

\heikodetail{

\bigskip

 {\tt MENTION SOMEWHERE IN THIS PAPER:
  The Function $F_{\textnormal{sym}}$ is unique with the property of having the same area as $F$.} see Lemma \ref{lem:equal-Cartan}.

  \bigskip

{\tt \yy Overath-Diss:
The $m$-harmonic symmetrization $F_{\sym}$ is the unique even and $1$-homogeneous function generating the same Finsler area functional as $F$. This is due to the representation of Finsler area in terms of the extended Radon transform as found in \eqref{cartan-radon} and the injectivity of the extended Radon transform on the set of smooth, $(-1)$-homogeneous and even real-valued functions. Notice, that the construction of Finsler area does not require the involved $F$ to be a full-fledged Finsler metric, but that $F$ is a (smooth), $1$-homogeneous function on $\R^n$, which satisfies $F(x,y)>0$ for all $x\in \R^n,\; y\in \R^n\backslash \{ 0 \}$ .}

\bigskip

}
\subsection{Main results}\label{sec:1.2}
{\bf Finsler-minimal graphs.}\,
As {\it Finsler-minimal graphs}  we denote  $m$-dimensional
Finsler-minimal immersions that
can be written as a graph over some domain in a hyperplane of the ambient 
space $\R^{m+1}$. Such a graph is called {\it entire} if the domain is the whole hyperplane.
\begin{theorem}[Bernstein theorems]\label{thm:bernstein}
Let $F=F(y)$ be a Minkowski metric on $\R^{m+1}$  satisfying
assumption {\rm (GA)}. Then the following holds:
\begin{enumerate}
\item[\rm (i)]
If $m=2$ or $m=3$, then every entire Finsler-minimal graph is an affine plane.
\item[\rm (ii)]
For every  $m\le 7$ there exists a constant $\delta=\delta(m)>0$ such that every
entire Finsler-minimal graph is an affine plane  if 
$$
\min\{\|F(\cdot)-|\cdot|\|_{C^3(\S^m)},
\|F_\sym(\cdot)-|\cdot|\|_{C^3(\S^m)}\}\le\delta.
$$
\item[\rm (iii)]
For every $m\in\N$ and $\gamma\in (0,1)$ there is a constant $\delta=\delta(m,\gamma) >0$
such that every entire Finsler-minimal graph
$\{u,f(u)):u\in\R^m\}\subset\R^{m+1}$ is an affine plane if 
\begin{equation}\label{C3-close}
\min\{\|F(\cdot)-|\cdot|\|_{C^3(\S^m)},
\|F_\sym(\cdot)-|\cdot|\|_{C^3(\S^m)}\}\le\delta,
\end{equation}
and if it satisfies the additional
growth condition
\begin{equation}\label{gradient-growth}
|Df(u)|=O(|(u,f(u)|^\gamma)\quad\textnormal{as $|u|\to\infty$}
\end{equation}
on the gradient $Df(u)=(f_{u^1},\ldots,f_{u^m}),$ where $u=(u^1,\ldots,u^m)\in\R^m$.
\end{enumerate}
\end{theorem}
Part (i) is apparently new for $m=3$, and it
generalizes \cite[Theorem 6]{sst} and \cite[Theorem 1.2]{cui-shen-2009} for
$m=2$ to general Finsler-Minkowski spaces. 
Indeed, calculating what (GA) implies for the specific Finsler metrics investigated
in \cite{sst} and \cite[Section 6]{cui-shen-2009}, we obtain exactly the threshold values
for the anisotropy $|b|$ mentioned above.
Part (ii) generalizes the part of \cite[Theorem 1.1]{cui-shen-2009} dealing with
Finsler-minimal graphs with respect to the Busemann-Hausdorff-area, since the given
scalar function $\phi(s):=(1+h(s))^{-1/m}$ with an arbitrary
odd smooth function $h:\R\to (-1,1)$ leads via an easy computation to an $(\alpha,\beta)$-Minkowski
metric satisfying (GA).
Part (iii) for arbitrary dimensions $m$ is, to the best of our knowledge,  
completely new in the context of Finsler-minimal immersions.

The next two results deal with uniqueness of Finsler-minimal graphs and the removability
of singularities.
\begin{theorem}[Uniqueness]\label{thm:uniqueness}
Let $F=F(y)$ be a Minkowski metric on $\R^{m+1}$ satisfying assumption {\rm (GA)}, and
assume that $f_1,f_2\in C^0(\bar{\Omega}\setminus K)\cap C^2(\Omega\setminus K)$ are
functions that define two Finsler-minimal graphs over $\Omega\setminus K\subset\R^m$,
where $\Omega$ is a bounded domain with a $C^1$-boundary $\partial\Omega$, and
$K\subset\Omega$ is compact with $\mathscr{H}^{m-1}(K)=0$ such that $\Omega\setminus
K$ is connected. Then, equality of $f_1$ and $f_2$ 
on the boundary $\partial\Omega$ implies 
$f_1=f_2$ on $\bar{\Omega}
\setminus K.$ 
\end{theorem}

The only known uniqueness result for Finsler-minimal graphs we are aware of  is
contained in the work of Souza et al.
\cite[Corollary 7]{sst} for the specific choice of a Randers-Minkowski 
$3$-space,
which can be recovered from Theorem \ref{thm:uniqueness} by choosing $m=2$
 and $K=\emptyset.$

 \begin{theorem}[Removability of singularities]\label{thm:removability}
 Let $F=F(y)$ be a Minkowski metric on $\R^{m+1}$ satisfying assumption {\rm (GA)}, 
 and assume that $\Omega\subset\R^m$ is an arbitrary domain and $K$
 is a locally compact subset of $\Omega$ with $\mathscr{H}^{m-1}(K)=0$.
 Then
 any Finsler-minimal graph of class $C^2$ on $\Omega\setminus K$
 can be extended as a Finsler-minimal graph of class $C^2$ onto all
 of $\Omega$.
 \end{theorem}
 With this result we generalize the only known removability result for
 Finsler-minimal graphs \cite[Proposition 12]{sst}, which can be reproduced
 from Theorem \ref{thm:removability} by setting $m=2$, $K:=\{u_0\}$ for some
 $u_0\in\Omega$, and choosing the particular Minkowski-Randers metric $F(y):=|y|+b_iy^i$
 with $|b|<1/\sqrt{3}.$ 
 
 Recall that our general assumption (GA) in Theorems 
 \ref{thm:bernstein} and \ref{thm:removability} is sharp since it translates
 to the sharp bounds on $|b|$ established in the work of Souza et al. and
 Cui and Shen in the specific Randers and $(\alpha,\beta)$-spaces. The minimal
  cones constructed in the respective spaces for $|b|$ above the threshold values
  demonstrate their sharpness; see \cite[p. 300]{sst} and \cite[Theorem 5.3]{cui-shen-2009}.

\bigskip

\noindent
{\bf Global results for Finsler-minimal immersions.}\,
Apart from Shen's nonexistence result for closed oriented Finsler-minimal submanifolds
in Minkowski space  \cite[Theorem 1.3]{shen98} and a few global results on rotationally
symmetric Finsler-minimal surfaces \cite{souza-tenenblat-2003}, \cite{cui-shen-2011}
in specific $(\alpha,\beta)$-spaces, there 
seems to be not much known about the global behaviour of Finsler-minimal immersions
that are not graphs.
We present here a simple enclosure theorem -- well known for classic
minimal surfaces in Euclidean space -- and a selection of isoperimetric inequalities
for Finsler-minimal immersions $X:\mathscr{M}\to\mathscr{N}:=\R^{m+1}$, where $\mathscr{M}$
denotes a smooth oriented $m$-dimensional manifold with 
boundary $\partial
\mathscr{M}$. 

\heikodetail{

\bigskip

\tt ORIENTED NECESSARY FOR ALL
FOLLOWING RESULTS?\xx\xx

Herr Overath bemerkt:
\yy Darueber muss ich nochmal nachdenken! Aber zumindest die meisten Resultate von Clarenz und Winklmann setzen orientiert oder orientierbar voraus, auch bei White wuerde ich vermute, dass dies sinnvoll ist. Allerdings genuegt es voraus zu setzen, dass die Ausgangsmannigfaltigkeit $\mathscr{M}$ orientiert ist.

\xx\xx Bei White gibt es aber auch eine Bemerkung, dass er $F(y)=F(-y)$ nur voraussetzt,
um keine Orientierung fordern zu muessen...\xx\xx

\yy\yy
Das ist richtig. Ich wuerde da allerdings gerne zumindest bei den anderen Resultaten auf Nummer sicher gehen und die Orientierbarkeit/Orientiertheit voraussetzen. Ggf. könnte zu dem White-Resultat eine entsprechende Bemerkung einfuegen!?!
\yy\yy
\bigskip

}
\begin{theorem}[Convex hull property]\label{thm:convex-hull}
Let $F=F(y)$ be a Minkowski metric on $\R^{m+1}$ such that assumption {\rm (GA)} holds.
Then the image $X({\mathscr{M}})$ of any Finsler-minimal immersion
$$
X\in C^2(\Int(\mathscr{M}),\R^{m+1})\cap C^0({\mathscr{M}},\R^{m+1})
$$
is contained in the convex
\heikodetail{

\bigskip

\tt Herr Overath bemerkt:
 \yy Checken?
Eine Anmerkung zum Begriff konvex: In Minkowskiraeumen ist dieser Begriff sinnvoll, da 
in diesen die Geodaeten ihrerseits wieder Strecken sind. In allgemeinen 
Finslerraeumen $(\R^{m+1}, F)$, also mit $x$-Abhaenigkeit, muss man zumindest ueber 
die Bedeutung des Wortes im Finslerschen nochmal nachdenken. Insbesondere sollten wir 
den Begriff definieren oder wenigstens etwas dazu anmerken. 
\xx Ist wohl jetzt mit dem Whiteschen Finsler-mean convex erledigt...
\bigskip

} 
hull of its boundary $X(\partial\mathscr{M})$, where
$\Int(\mathscr{M})=\mathscr{M}\setminus \partial\mathscr{M}.$
\end{theorem}

For a general Finsler metric $F=F(x,y)$ we set
\begin{equation}\label{definite}
M_F:=\sup_{\R^{m+1}\times\S^m}F(\cdot,\cdot)\AND m_F
:=\inf_{\R^{m+1}\times\S^m}F(\cdot,\cdot).
\end{equation}
Notice that for a Minkowski metric one automatically has $0<m_F\le M_F<\infty$, since
(F1) and (F2) imply that $F>0$.
A simple variant of the isoperimetric inequality  of
 the form 
$$
\area^F(X)\le  \frac{M_F^2}{4\pi m_F^2}\Big(\mathscr{L}^F(\Gamma)\Big)^2
$$
has been shown for (possibly branched) Finsler-area {\it minimizing} surfaces in $\R^3$
with a given boundary contour $\Gamma\subset\R^3$ by a simple comparison
with classic minimal surfaces; see  \cite[Corollary 1.3]{overath-vdM-2012a}. Here,
$\mathscr{L}^F(\Gamma):=\int F(\Gamma,\dot{\Gamma})$ denotes the Finslerian length
of $\Gamma$. 
 
For the following isoperimetric inequalities for Finsler-minimal immersions
in Finsler space denote by $dS_F$ the volume form on the
boundary $\partial\mathscr{M}$
induced by $X$ whose restriction $X|_{\partial\mathscr{M}}:\partial\mathscr{M}\to\R^{m+1}$
is again an immersion. Analogous to \eqref{volumefactor} one defines
$$
\int_\omega dS_F:=\area_\omega^F(X|_{\partial\mathscr{M}})
$$
for any relatively open set $\omega\subset\partial\mathscr{M}.$
For the special choice $F(\cdot)=|\cdot|$, 
we set $\int_{\Omega}\,dS:=\int_{\Omega}\,dS_F$.
\begin{theorem}[Isoperimetric inequalities]\label{thm:isop}
Let $F=F(x,y)$ be a Finsler metric on $\R^{m+1}$ satisfying assumption {\rm (GA)}
and let $X\in C^2(\Int(\mathscr{M}),\R^{m+1})
\cap C^1(\mathscr{M},\R^{m+1})$ be a Finsler-minimal immersion.
Then the following holds.
\begin{enumerate}
\item[\rm (i)] If $F=F(y)$ and 
 $X(\partial\mathscr{M})\subset \overline{
B_R^F(a)}$ for some $a\in\R^{m+1}$, where $\overline{B^F_R(a)}$ denotes the closed Finsler-Minkowski unit ball $\lbrace v\in\mathbb{R}^{m+1}: 
F(v-a) \le R\rbrace$,
then 
\begin{equation}\label{isop1}
\area^F_\mathscr{M}(X)\le\frac{R}{m}\left(\frac{M_F}{m_F}\right)^m\sqrt{1+
\Lambda(F)(m+1)(m/m_F)^2}\int_{\partial\mathscr{M}}\,dS_F,
\end{equation}
where $\Lambda(F)$ denotes the largest possible eigenvalue of the fundamental tensor
$(g_{ij})$ of $F$ when restricted to the sphere $\S^m$.
\item[\rm (ii)]
If $m=2$ and $F=F(y)$, and
if the boundary $\partial\mathscr{M}$ consists of $k\ge 1$ closed
rectifiable Jordan curves $\gamma_i$ with images $\Gamma_i:=X(\gamma_i)$ for $
i=1,\ldots,k$, then for all $a\in\R^3$ 
\begin{equation}\label{isop2}
\area^F_\mathscr{M}(X)\le \frac{M_F^2}{m_F^2}\sqrt{1+\Lambda(F)\frac{12}{m_F^2}}
\sum_{i=1}^k\Big[\frac{\mathscr{L}^F(\Gamma_i)^2}{4\pi}
+\frac 12 \mathscr{L}^F(\Gamma_i)\dist_F(a,\Gamma_i)\Big],
\end{equation}
where $\dist_F$ denotes the Finslerian distance in $\R^3$ induced by $F$, and
where 
$\Lambda(F)$ is as in part {\rm (i).}
\item[\rm (iii)]
Let $m=2$, $\mathscr{M}=B\equiv B_1(0)\subset\R^2$, $X\in C^{2,\alpha}(\bar{B},\R^3)$
with $\|X\|_{L^\infty(B,\R^3)}\le 1$,
mapping $\partial B$ topologically onto the closed Jordan curve $\Gamma\subset
B_R(a)\subset\R^3$ for some $a\in\R^3$, $R\le 1$. Then there is a universal constant $\delta>0$ such that
\begin{equation}\label{isop3}
\area^F_B(X)\le R\left\{c_1(F){M_F^*}^2\Big[\int_\Gamma\kappa\,ds-2\pi\Big]+
c_2(F)\frac{{M_F^*}^2}{m_F^*}\mathscr{L}^F(\Gamma)\right\},
\end{equation}
as long as 
\begin{equation}\label{C3-Finsler-close}
\min\{\|F(\cdot)-|\cdot|\|_{C^3(\overline{B_1(0)}\times\S^m)},
\|F_\sym(\cdot)-|\cdot|\|_{C^3(\overline{B_1(0)}\times \S^m)}\}\le\delta
\end{equation}
holds. Here $M_F^*$ and $m_F^*$ denote  the supremum and the  infimum of $F$
on $\overline{B_1(0)}\times \S^m$, respectively, and $\int_\Gamma\kappa\,ds$ is the
Euclidean total curvature of the boundary curve $\Gamma$.
\end{enumerate}
\end{theorem}
Notice for parts (i) and (ii) that $M_F=m_F=\Lambda(F)=1$ 
if $F$ equals the Euclidean norm, that is, $F(x,y)=|y|$.
Likewise in part (iii) one can show, that the constant $c_1(F)$ vanishes and
$c_2(F)=1/2$, and $M_F^*=m_F^*=1$, if $F$ is the Euclidean norm.
On the other hand, if the boundary contour $\Gamma$ happens to be planar and
convex, then according to Fenchel's theorem the first summand vanishes
 in \eqref{isop3} and one is left with a fairly simple isoperimetric
 inequality as long as $F$ satisfies \eqref{C3-Finsler-close}.

\bigskip

\noindent
{\bf Existence and uniqueness of Finsler-minimal immersions spanning
given boundary contours.}\,
For a Finsler metric $F=F(x,y)$ on $\R^3$ satisfying our general assumption (GA)
and for a given rectifiable Jordan curve
$\Gamma\subset\R^3$, we could establish  in \cite[Theorem 1.2]{overath-vdM-2012a} 
the existence
of conformally parametrized minimizers of $\area^F_B$ in the class  of Sobolev
mappings $W^{1,2}(B,\R^3)$ from the two-dimensional unit ball $B=B_1(0)\subset\R^2$
into $\R^3$ that parametrize $\Gamma$ on $\partial B$ in a weakly monotonic way\footnote{See
\cite[pp. 231--232]{DHKW1} for the notion of weakly monotonic mappings 
on the boundary.}. The proof produces a mild improvement in regularity of
the $\area^F_B$-minimizer, i.e., continuity
up to the boundary,  interior H\"older continuity and a slightly better integrability
of the gradient, but in general there is up to now no way to exclude branch points,
or even to estimate the size of the set of branch points. These are points, where
the Jacobian of the mapping fails to have rank two. If one asks for immersed
 or even embedded 
Finsler-minimal surfaces spanning a given boundary contour one has to add more
assumptions on the boundary curve.
\begin{theorem}[Existence]\label{thm:existence}
Let $F=F(x,y)$ be a Finsler metric on $\R^3$ satisfying assumption {\rm (GA)}. 
\begin{enumerate}
\item[\rm (i)]
If $\Gamma\subset\R^3$ is a smooth closed Jordan curve contained in the boundary
of a strictly Finsler area mean convex body, then for each $g\ge 0$ there exists a smooth embedded Finsler-minimal surface
spanning $\Gamma$ with genus less or equal $g$. 
In particular, such an embedded Finsler-minimal surface exists 
if $F=F(y)$ and if $\Gamma$ is contained in the boundary of a 
strictly convex body. 
\item[\rm (ii)]
If $F=F(y)$ and $\Gamma$ is a graph of bounded slope over $\partial\Omega$ for some bounded
convex domain $\Omega\subset\R^2$, then there is a smooth
and (up to reparametrizations)
 unique Finsler-minimal graph spanning $\Gamma.$
\end{enumerate}
\end{theorem}

Roughly speaking, a body is Finsler mean convex 
if inward variations of its boundary 
lead to an infinitesimal decrease of Finsler area; for details we refer to
Section \ref{sec:3}.
The bounded slope condition in the last part of the theorem means that we
find a constant $R>0$, so that we 
can write $\Gamma$ as a graph,
$$
\Gamma=\{(u,\gamma(u))\in\R^3: u=(u^1,u^2)\in\partial\Omega\}
$$
for some function $\gamma:\partial\Omega\to\R$, 
such that for any curve point
$(u_0,\gamma(u_0))\in\Gamma$  there exist
two vectors $p^+_0,p^-_0\in \overline{B_R(0)}\subset\R^2$
\heikodetail{

\bigskip

Herr Overath bemerkt:
\tt \yy Hier koennte man $\overline{B_R(0)}$ durch $\overline{B^F_R(0)}$ ersetzen, 
um die Bedingung Finslersch zu formulieren.

Aber davon haben wir abgesehen..

\bigskip

} 
such
that the two affine linear functions 
$$
\ell^+_0(u):= p^+_0\cdot (u-u_0)+\gamma(u_0)\AND
\ell^-_0(u):= p^-_0\cdot (u-u_0)+\gamma(u_0)
$$
satisfy
$\ell^-_0(u)\le\gamma(u)\le\ell^+_0(u)$ for all $u\in\partial\Omega$. 
In particular, if $\Omega$ is strictly convex and $\Gamma$ is a $C^2$-graph over 
$\partial\Omega$ then $\Gamma$ satisfies the bounded slope  condition; see
\cite[pp. 309, 310]{gilbarg-trudinger_1998}. 

For part (i) of Theorem \ref{thm:existence} one can even prescribe an upper
bound $G\ge 0$ on the genus of the Finsler-minimal embedding, which itself
is Finsler-area minimizing among all embedded surfaces with the same boundary
and with genus less or equal $G$.

\bigskip

\noindent
{\bf Strategy of proofs and structure of the paper.}
The key to proving all these results is a connection between Finsler minimal surfaces
and so-called Cartan functionals, i.e., parameter invariant variational integrals
with a specific structure of the integrand. This connection has been established
in \cite{overath-vdM-2012a} by means of the spherical Radon transform. 
We will recall all relevant facts on Cartan functionals and on
their relation to Finsler minimal immersions in the next section; see Section
\ref{sec:2.1}. In addition, in order to apply the existing theory on
Cartan functionals 
we need to strengthen this connection to
show that Finsler area generates an {\it elliptic} Cartan integrands, that is, a strict
parametric convexity property that was not necessary for the weak existence
theory for the Plateau problem presented in \cite{overath-vdM-2012a}. Here, however,
we need to analyse the behaviour of the spherical Radon transform on
a suitable function space endowed with a  Fr\'echet topology; see Section \ref{sec:2.2},
in particular Corollary \ref{cor:ellipticity} ensuring elliptic Cartan integrands
related to Finsler area. In Section  \ref{sec:3} we apply various results on
critical immersions of Cartan functionals established by 
H. Jenkins \cite{jenkins_1961b},
L. Simon \cite{simon_1977a}, \cite{simon_1977b}, B. White \cite{white-existence-1991},
U. Clarenz and the second author \cite{clarenz_enclosure}, 
\cite{clarenz-vdM-compact_2001}, \cite{clarenz-vdM-isop},
S. Winklmann \cite{winklmann-isop}, \cite{winklmann-bernstein},
\cite{winklmann-existence-uniqueness}, and S. Hildebrandt and F. Sauvigny,
\cite{Hil-sauvigny-energy-estimate-2009} to prove the results stated above.

{\bf Acknowledgments.} Substantial parts of this work are contained in the first author's thesis who was partially supported by DFG grant no. Mo 966/3-1,2, and moreover
by the Excellence Initiative of the German federal and state governments.
Moreover, the second author would like to express his gratitude to Professor
Seiki Nishikawa for inviting him to Tohoku University at Sendai, Japan, to
give a series of talks on this line of research.

\setnumbers
\section{Finsler-minimal immersions from a variational viewpoint}
\label{sec:2}
\subsection{Finsler area and Cartan functionals}\label{sec:2.1}

The explicit form \eqref{volumefactor} of the Busemann-Hausdorff volume implies that
bounds on the Finsler metric directly transfer to corresponding bounds on the area functional.
\begin{lemma}[Area comparison]\label{lem:area-comparison}
 Let $F=F(x,y)$ be a Finsler metric on $\R^n$, and denote by $E=E(y):=|y|$ the Euclidean metric on $\R^n$. Assume that there are two constants $0<m_F\le M_F$ such that
\begin{equation}
 m_F|y|=m_FE(y) \le F(x,y) \le M_FE(y)=M_F|y|\label{lem:area-comparison-eq1}
\quad\Foa (x,y)\in \R^n\times \R^n.
\end{equation}
 Then one has
\begin{equation}
 m_F^k \area^{E}_{\Omega}(X) \le \area^{F}_{\Omega}(X) \le M_F^k \area^{E}_{\Omega}(X)\label{lem:area-comparison-eq2}
\end{equation}
for every immersion $X:{\Sigma^k}\rightarrow \R^n$ from a smooth 
$k$-dimensional manifold ${\Sigma}$ into $\R^n$ and for every open set $\Omega\subset\Sigma$.
\end{lemma}
\proof
Inequality \eqref{lem:area-comparison-eq1} implies $m_FX^*E(v)\le X^*F(u,v)\le 1$ for all $u\in\Omega$ and all $v\in {\lbrace v=(v^1,\ldots,v^k)\in T_u\Omega\simeq \R^k : X^*F(u,v)\le 1\rbrace}$. Hence we get
\begin{equation}
{\lbrace v\in \R^k : X^*F(u,v)\le 1\rbrace} \subset {\lbrace v\in \R^k : X^*E(u,v)\le \frac{1}{m_F}\rbrace}.\label{proof:area-comparison-Finsler-unit-ball-inclusion}
\end{equation}
The $k$-dimensional Hausdorff-measure is monotonic 
\heikodetail{

\bigskip

two lines above we have used in detail:
$$
m_FX^*E(u,v)=m_FE(X(u),dX|_u(v))\overset{\eqref{lem:area-comparison-eq1}}{\le}
F(X(u),dX|_u(v))=X^*F(u,v)
$$

monotonic means that for $A\subset B$ one has $\H^k(A)\le \H^k(B)$

and for the following homogeneity the details read as:

\begin{eqnarray*}
X^*(E(u,tv) & = &
E(X(u),dX|_u(tv)) = E(X(u),tdX|_u(v))\\
& \overset{\eqref{H}}{=} & tE(X(u),dX|_u(v))=tX^*(E(u,v).
\end{eqnarray*}

\bigskip

}
and scales like 
$
 \mathscr{H}^k(tA) = t^k\mathscr{H}^k(A)
$
for all $t>0$ and for any measurable set $A\subset \R^n$; see, e.g., \cite[Chapter 2.1, Theorem 2]{evans-gariepy}. Therefore \eqref{proof:area-comparison-Finsler-unit-ball-inclusion} yields
\begin{eqnarray}
\mathscr{H}^k({\lbrace v\in \R^k : X^*F(u,v)\le 1\rbrace})
& \le  &\mathscr{H}^k({\lbrace v\in \R^k : X^*E(u,v)\le \frac{1}{m_F}\rbrace})\nonumber\\
& = & \mathscr{H}^k(\frac{1}{m_F}{\lbrace v\in \R^k : X^*E(u,v)\le 1\rbrace})\nonumber\\
& = & \frac{1}{m_F^k}\mathscr{H}^k({\lbrace v\in \R^k : X^*E(u,v)\le 1\rbrace}).\label{proof:area-comparison-Finsler-unit-ball-measure-estimate}
\end{eqnarray}
Combining \eqref{proof:area-comparison-Finsler-unit-ball-measure-estimate} with \eqref{volumefactor} leads to
$
 m_F^k \sigma_{X^*E}(u) \le \sigma_{X^*F}(u)$,
and by means of the definition \eqref{area} of Finsler area we thus obtain
the left inequality in \eqref{lem:area-comparison-eq2}. The second estimate in \eqref{lem:area-comparison-eq2} can be shown analogously by exchanging the role of $E$ and $F$ and exploiting the second estimate in \eqref{lem:area-comparison-eq1}.

\qed

\bigskip

Recall from \cite[Theorem 1.1]{overath-vdM-2012a} that 
we can rewrite Finsler area for $\mathscr{N}=\R^{m+1}$ as follows.
\begin{theorem}\label{thm:finsler-area}
If $\mathscr{N}=\R^{m+1}$ with a Finsler structure $F=F(x,y)$, then
the Finsler area of any $C^1$-immersion $X$ of a smooth $m$-dimensional
manifold $\mathscr{M}$ into $\R^{m+1}$ may be expressed in local
coordinates $(u^1,\ldots,u^m):\Omega\subset\mathscr{M}\to
\tilde{\Omega}\subset\R^m$ as
\begin{equation}\label{area_space}
\area^F_\Omega(X)=\int_{\tilde{\Omega}}\cA^F(X(u),\big(\frac{\partial X}{\partial u^1}
\wedge
\ldots\wedge\frac{\partial X}{\partial u^m}\big)(u))\,du^1\wedge\ldots\wedge du^m,
\end{equation}
where  
\begin{equation}\label{AF_space}
\cA^F(x,Z)=\frac{|Z|\mathscr{H}^m(B_1^m(0))}{\mathscr{H}^m(\{T\in Z^\perp\subset\R^{m+1}:
F(x,T)\le 1\})}\quad\Fo (x,Z)\in\R^{m+1}\times (\R^{m+1}\setminus\{0\}). 
\end{equation}
\end{theorem}

Notice that  the explicit form \eqref{AF_space}  of the integrand $\cA^F$
directly implies positive $1$-homogeneity in its second argument: 
\begin{equation}\label{H}
\cA^F(x,tZ)=t\cA^F(x,Z)\quad\Foa (x,Z)\in\R^{m+1}\times(\R^{m+1}\setminus\{0\}), \,t>0.
\end{equation}
This together with the fact that $\cA^F$ in the integral in \eqref{area_space}
depends
on the position $X(u)$ and
on the normal vector 
$\left(X_{u^1}\wedge\ldots\wedge X_{u^m}\right)(u)$  qualifies
$\cA^F$ as {\it Cartan integrand} as defined in \cite[p. 2]{HilvdM-parma}. 
Moreover, inspecting \eqref{AF_space} one  immediately sees that $\cA^F$
does not depend on $x$ in case of a Minkowski metric $F=F(y)$, and
$\cA^F(x,Z)$ simplifies to the Euclidean area integrand $\cA(Z)=|Z|$
if $F$ happens to be the Euclidean metric, $F(x,y)=E(y)=|y|$. Notice, in addition,
that the Cartan integrand $\cA^F$ is even in its second argument even if
$F$ is not reversible. Finally, we will see in Corollary \ref{cor:radon-continuous}
below that $\cA^F$ is smooth on $\R^{m+1}\times (\R^{m+1}\setminus\{0\})$
since $F$ is smooth on this set. 

\heikodetail{

\bigskip
IST NUN DRIN:

{\tt\xx\xx\xx ALL DAS FOLGENDE MUSS NOCH PLATZIERT UND AUSGEFUEHRT WERDEN!!}
In addition, \yy$\cA^F(x,Z)$ is an even function in the $Z$-variable, \yy even when $F$ is not reversible \yy.{\tt \yy $\cA^F$ is always even, a result of the fact, that the Radon transform annihilates odd functions.\yy}
{\tt\xx VIELLEICHT DOCH BESSER, DIE EXPLIZITE FORM ALS QUOTIENT DER MASSE ANZUGEBEN\xx}{\tt \yy Koennte man sicher machen, auch im Verweis auf die vorige Arbeit, andererseits stehen die Definitionen der Ausdruecke im wesentlichen in der Einleitung, so dass man vielleicht einfach die Ausdruecke hier, damit in Beziehung setzen koennte.} 
{\tt\xx BAUE AUCH NOCH DIE GLATTHEIT VON $\cA^F$ EIN, DIE AUS GLATTHEIT VON $F$
FOLGEN SOLLTE...\xx}{\tt\yy Die Glattheit kann man aus der Darstellung mittels der Radon transform herausziehen, welche wir in \cite{overath-vdM-2012a} diskutiert haben.}

\bigskip

}

Even if $F$ itself is not a Finsler metric one can deduce several useful 
properties of $\cA^F$, for example, explicit $L^\infty$-bounds on  $F=F(x,y)$
transfer to corresponding bounds on the Cartan area integrand $\cA^F=\cA(x,Z)$
proven in \cite[Lemma 2.4]{overath-vdM-2012a}:
\begin{lemma}[Pointwise bounds]\label{lem:pointwise-bounds}
Let $F_1(x,y),F_2(x,y)$ be continuous on $\R^{m+1}\times\R^{m+1}$,
strictly positive for $y\not=0$, and positively $1$-homogeneous in the $y$-variable,
and assume that for each $x\in\R^{m+1}$ there exist numbers $0<c_1(x)\le
c_2(x)$ such that
\begin{equation}\label{F-bounds}
c_1(x)F_1(x,y)\le F_2(x,y)\le c_2(x)F_1(x,y)\quad\Foa y\in\R^{m+1},
\end{equation}
then the corresponding Cartan area integrands $\cA^{F_1}$ and $\cA^{F_2}$
satisfy
\begin{equation}\label{A-bounds}
c_1^m(x)\cA^{F_1}(x,Z)\le\cA^{F_2}(x,Z)\le c_2^m(x)\cA^{F_1}(x,Z)\quad\Foa Z\in
\R^{m+1}.
\end{equation}
\end{lemma}
If one wants to use this bridge between Finsler area and the variational theory
for Cartan functionals to prove new results about Finsler-minimal immersions, one 
has to establish convexity of $\cA^F$ in the $Z$-variable. For reversible Finsler
metrics this was indeed proven geometrically by Busemann 
\cite[Theorem II, p. 28]{busemann-convex-brunn-minkowski-1949}
in the completely different context of volume computations for cross sections
of convex bodies:
\begin{theorem}[Busemann]\label{thm:busemann}
If $F$ is a reversible Finsler metric on $\mathscr{N}=\R^{m+1}$, then the corresponding
Cartan area integrand $\cA^F=\cA^F(x,Z)$ is convex in the $Z$-variable for
any $x\in\R^{m+1}.$
\end{theorem}
It was Busemann's essential requirement of {\it reversible} Finsler metrics that motivated
our choice of symmetrization for in general non-reversible Finsler structures
$F$. That the $m$-harmonic symmetrization defined in \eqref{m-harmonic}
is indeed suitable is partially justified by the following result
shown in \cite[Lemma 2.3]{overath-vdM-2012a}:
\begin{lemma}\label{lem:equal-Cartan}
If $F=F(x,y)$ is continuous on $\R^{m+1}\times\R^{m+1}$, strictly
positive whenever $y\not=0$, and positively homogeneous in the $y$-variable,
then
\begin{equation}\label{equal-Cartan}
\cA^F(x,Z)=\cA^{F_\sym}(x,Z)\quad\Foa (x,Z)\in\R^{m+1}\times\R^{m+1}.
\end{equation}
\end{lemma}
Let us mention here that the $m$-harmonic symmetrization $F_\sym$
is the
unique even and positively  
$1$-homogeneous function generating the same area integrand
as in \eqref{equal-Cartan}. This is due to the invertibility of the
extended spherical Radon transform  on  even
positively $(-m)$-homogeneous smooth functions on $\R^{m+1}\setminus\{0\}$;
see Theorem \ref{thm:inverse-radon} and Lemma \ref{lem:cartan-radon}
in Section \ref{sec:2.2}. 

Convexity together with the pointwise bounds stated in Lemma \ref{lem:pointwise-bounds} suffices to apply the existence theory for Cartan functionals, which helped
solving the Plateau problem in Finsler $3$-space in \cite{overath-vdM-2012a}.
But in order to access
available results on critical immersions for Cartan functionals one needs
strict convexity of $\cA^F$ in the $Z$-variable in the sense of {\it parametric ellipticity}. This is the 
strongest from of convexity that one can expect due to the homogeneity
in the $Z$-variable, as described in the next subsection. Let us point out
that we do not see how to quantify Busemann's original geometric proof to obtain
strict convexity of the Cartan area integrand, which is why we devised an alternative route via the spherical
Radon transform.

\subsection{Radon transform and parametric ellipticity}\label{sec:2.2}
The {\it spherical Radon transform}\footnote{For $m=2$ also known as the
{\it Funk transform} \cite{funk-1915}.} \cite{radon-1994}
functions that are continuous on the unit sphere:
\begin{definition}[Spherical Radon transform]\label{def:radon-transform}
The {\em spherical Radon transform} $\widehat{\cR}[f]$ of a function
$f\in C^0(\S^m)$ is defined as
\begin{equation}\label{spherical-radon}
\widehat{\cR}[f](\zeta):=
\frac{1}{\H^{m-1}(\S^{m-1})}\int_{\S^m\cap\zeta^\perp}f(\omega)\,d\H^{m-1}(\omega)
\quad\Fo \zeta\in\S^m.
\end{equation}
\end{definition}
In the context of integral geometry, geometric tomography, and convex analysis
the spherical Radon transform has been used intensively; see, e.g., 
\cite{helgason_radon, helgason_groups,bailey,gardner,groemer}. In our present
setting it is useful to look at a suitable homogeneous extension of the Radon
transform. We define
\begin{equation}\label{R}
\cR[g](Z):=\frac{1}{|Z|}\hR\big[g|_{\S^m}\big]\Big(\frac{Z}{|Z|}\Big)
\quad\Fo
g\in C^0(\R^{m+1}\setminus\{0\}),\,\,\, Z\in\R^{m+1}\setminus\{0\},
\end{equation}
which by definition is a $(-1)$-homogeneous function on $\R^{m+1}\setminus\{0\}$.

In \cite[Section 3]{overath-vdM-2012a} we analyzed in detail this
homogeneous extension and proved among other things that $\cR$ is a bounded
linear map from $C^0(\R^{m+1}\setminus\{0\})$ to itself \cite[Corollary 3.3]{overath-vdM-2012a}, and that $\cR[g]$ 
is of class $C^k$ if $g$ is $(-m)$-homogeneous and itself of class $C^k$, together
with an explicit differentiation formula \cite[Theorem 3.6]{overath-vdM-2012a},
to produce quantitative sufficient criteria to guarantee higher regularity for
the solutions of the Plateau problem in Finsler $3$-space. 
Here, we will need this differentiation rule only
once, and only for the first order derivative:
\begin{equation}\label{diff-rule}
Z_\tau\frac{\partial}{\partial Z_\sigma}\cR[g](Z)=-\cR\Big[
\frac{\partial}{\partial y^\tau}(y^\sigma g)\Big](Z)
\end{equation}
for all $Z=(Z^1,\ldots,Z^{m+1}),$\, $y=(y^1,\ldots,y^{m+1})\in\R^{m+1}
\setminus \{0\}$, where we set $Z_{\tau}:=\delta_{\tau l}Z^l.$

The key observation
to connect the Radon transformation to Finsler-minimal surfaces is,
that one can rewrite the Cartan area integrand $\cA^F$ in terms of
the extended spherical Radon transform; see \cite[Corollary 3.8]{overath-vdM-2012a}:
\begin{lemma}[Cartan area as Radon transform]
\label{lem:cartan-radon}
For any function $F=F(x,y)$ that is continuous on $\R^{m+1}\times\R^{m+1}$,
positive whenever $y\not=0$, and positively homogeneous in the $y$-variable,
one has the identity
\begin{equation}\label{cartan-radon}
\cA^F(x,Z)=\frac{1}{\cR\big[F^{-m}(x,\cdot)\big](Z)}\quad\Fo
(x,Z)\in\R^{m+1}\times (\R^{m+1}\setminus\{0\}).
\end{equation}
\end{lemma}
Our main goal in this section is to prove that {\it every} Finsler metric $F$ leads to
a Cartan area integrand $\cA^F(x,Z)$ that is strongly convex on $Z^\perp$
in the second variable, which corresponds to {\it parametric ellipticity};
see Definition \ref{def:ellipticity} below.
In order to do
this we are going to use an indirect reasoning, which requires to invert
 the Radon transform to go back and forth between the original Finsler
metric  $F$ and its image $\cA^F$ under the Radon transform. Such an 
invertibility result is available only on the Fr\'echet space of even
smooth functions
on the sphere, since one gains some orders of derivatives under the Radon
transform. In \cite[Theorem 3.4.14 \& Proposition 3.6.4]{groemer} 
one finds a quantitative
estimate on the inverse proven by means of spherical harmonics that can be extended to 
all orders of differentiability; see \cite[Prop. 2.2.26 \& Theorem 2.2.27]{overath-phd_2013}.
\heikodetail{

\bigskip

Herr Overath schreibt:
\tt\yy\yy Die Verallgemeinerung von Satz \cite[Theorem 3.4.14]{groemer} auf alle Ableitungsstufen wurde in \cite{overath-phd_2013} nicht direkt fuer die sphaerische 
Radon trafo sondern fuer die erweiterte gezeigt, es folgt aber 
aus  \cite[Proposition 2.2.26 \& Theorem 2.2.27]{overath-phd_2013} 
und den in deren Beweisen enthaltenen Abschaetzungen durch 
Einschraenkung auf die Sphaere.\yy

\bigskip

}
But  for our purposes it will suffice to prove continuity of the extended
Radon transform $\cR$ on a certain Fr\'echet subspace of
$C^\infty(\R^{m+1}\setminus\{0\})$,
and to give an explicit formula of the inverse $\cT:=\cR^{-1}$ in terms
of Helgason's inverse of the spherical Radon transform 
$\widehat{\cT}={\widehat{\cR}}^{-1}$ as, e.g., 
presented in \cite[Chapter III, Theorem 1.11]{helgason_radon}. 
Combining this with the open mapping theorem on Fr\'echet
spaces yields the desired (non-quantitative) continuity of the inverse $\cT.$

To begin with, recall from \cite[Corollary 3.7]{overath-vdM-2012a} the following
estimate of seminorms
\begin{equation}\label{semi-norms-sphere}
\varrho_l(g):=\max\{|D^\alpha g(\xi)|:\xi\in\S^m,|\alpha|\le l\}\quad\Fo
l=0,1,\ldots,k,
\end{equation}
on functions $g$ of class $C^k(\R^{m+1}\setminus\{0\})$, where $\alpha=(\alpha_1,\ldots,
\alpha_m)$ is
a multi-index with $\alpha_i\in\N\cup\{0\}$ and with
length $|\alpha|=\alpha_1+\cdots +\alpha_m$:
\begin{lemma}\label{lem:semi-norms-sphere}
There is a constant $C=C(m,k)$ such that for any positively $(-m)$-homogeneous
function $g\in C^k(\R^{m+1}\setminus\{0\})$ one has
\begin{equation}\label{semi-norms-sphere-estimate}
\varrho_k(\cR[g])\le C(m,k)\varrho_k(g).
\end{equation}
\end{lemma}
We can use \eqref{semi-norms-sphere-estimate} to estimate the seminorms 
\begin{equation}\label{semi-norms-space}
p_{k}(f):=\max\{|D^\alpha f(Z)|: Z\in S_k,\,|\alpha|\le k\},\quad k=0,1,2,\ldots,
\end{equation}
on the nested compact sets
$$
S_k:=\overline{B_{k+1}(0)}\setminus B_{\frac{1}{k+1}}(0)\subset\R^{m+1}
$$
for positively homogeneous functions:
\begin{lemma}\label{lem:relation-semi-norms}
For a positively $q$-homogeneous function $f\in C^k(\R^{m+1}\setminus\{0\})$ one
has
\begin{equation}\label{relation-semi-norms}
\varrho_k(f)\le p_{k}(f)\le (k+1)^{|q|+k}\varrho_k(f)\quad\textnormal{for all $k=0,1,2,\ldots$.}
\end{equation}
\heikodetail{

\bigskip

Herr Overath bemerkt:
\tt \yy Anmerkung: Dafuer habe ich in meiner Doktorarbeit $2$er Potenzen gewaehlt, so wie in Rudin, denke ich.

\bigskip

}
\end{lemma}
\proof
The first estimate is trivial since $\S^m\subset S_k$ for all $k=0,1,2,\ldots .$
Notice for the second inequality that $D^\alpha f$ is positively
$(q-|\alpha|)$-homogeneous, i.e.,
$$
D^\alpha f(Z)=|Z|^{q-|\alpha|}D^\alpha f(Z/|Z|)\quad\Fo Z\in\R^{m+1}\setminus\{0\},
$$
 so that
 \begin{eqnarray*}
 |D^\alpha f(Z)| & \le & |Z|^{q-|\alpha|}\max\{|D^\alpha f(\xi)|:
 \xi\in\S^m\}\\
 & \le & |Z|^{q-|\alpha|}\varrho_k(f)\quad\Fo |\alpha|\le k.
 \end{eqnarray*}
 For $Z\in S_k$ with $|Z|\ge 1$ one has 
 $
 |Z|^{q-|\alpha|}\le (k+1)^{|q|+|\alpha|},$ and for $1/(k+1)\le |Z|<1$ the estimate
 $$
 |Z|^{q-|\alpha|}\le|Z|^q(k+1)^{|\alpha|}=\begin{cases} {(k+1)^{|\alpha|}}{|Z|^{-|q|}} &
 \Fo q<0,\\
 |Z|^{|q|}(k+1)^{|\alpha|} & \Fo q>0,
 \end{cases}
 $$
 and the right-hand side is clearly dominated by $(k+1)^{|q|+|\alpha|}$ as well.
 This leads to
 $$
 |D^{\alpha}f(Z)|\le (k+1)^{|q|+k}\varrho_k(f)\quad\Foa Z\in S_k,\,
 |\alpha|\le k,
 $$
 for all  $k=0,1,2,\ldots,$ which concludes the proof.
 \qed

 Since $\R^{m+1}\setminus\{0\}=\bigcup_{k=1}^\infty S_k$, the space 
 $C^\infty(\R^{m+1}\setminus\{0\})$ equipped with the family of seminorms
 $p_k$ for $k=0,1,2,\ldots,$  is a Fr\'echet space; see 
 \cite[Section 1.46]{rudin-fa-book_1973}. The closed subspace of functions
 $g\in C^\infty(\R^{m+1}\setminus\{0\})$ that are, in addition, positively
 $(-m)$-homogeneous is therefore itself a Fr\'echet space with respect to
 the same family of seminorms, 
 \heikodetail{
 
 \bigskip
 
 \tt\xx eine nicht-optimale Quelle daf\"ur w\"are 
 \cite[Section 18, 3.]{koethe-buch}, da bemerkt er, dass Teilr\"aume von lokalkonvexen
 R\"aumen wieder lokalkonvex sind. Zusammen mit der Tatsache, dass abgeschlossene
 Unterr\"aume von vollst\"andigen metrischen R\"aumen wieder vollst\"andig sind,
 h\"atte man dann auch wohl dieses statement, im Rudin selbst habe ich kein
 statement gefunden, dass direkter zitierbar w\"are...wenn wir \"uberhaupt
 hier etwas zitieren wollen...?\xx 
 
 \yy Ihre Argumentation ueber \cite[Section 18, 3.]{koethe-buch} erscheint mir korrekt. Im Rudin wird dies nur indirekt angemerkt in \cite[Example 1.45]{rudin-fa-book_1973}. Im englischen Wikipedia-Artikel zu Frechet spaces steht sogar explizit drin, dass abgeschlossene Unterraeume von Frechet spaces ihrerseits Frechet spaces sind. Ggf. kann ich die zugehoerigen Quellen nochmal durchgehen, insofern Sie die Aussage durch unsere Quellen als nicht hinreichend belegt sehen. 
 
 \bigskip
 
 } 
and we can
 combine
the preceding two lemmas to obtain the continuity of the extended Radon
 transform on  this smaller Fr\'echet space.
 \begin{corollary}[Continuity of $\cR$]\label{cor:radon-continuous}
 The extended Radon transform $\cR$ is a bounded linear mapping
 from the Fr\'echet space of $(-m)$-homogeneous mappings of class
 $C^\infty(\R^{m+1}\setminus\{0\})$ into the Fr\'echet space
 of $(-1)$-homogenous mappings of class
 $C^\infty(\R^{m+1}\setminus\{0\})$, satisfying the estimate
 \begin{equation}\label{radon-continuous}
 p_k(\cR[f])\le C(m,k)(k+1)^{k+1}p_k(f)\quad k=0,1,2,\ldots,
 \end{equation}
 where $f\in C^\infty(\R^{m+1}\setminus\{0\})$ is an arbitrary
 positively $(-m)$-homogeneous function.
 \end{corollary}
 \proof
 That $\cR$ is linear and that $\cR[f]$ is positively $(-1)$-homogeneous
 can be seen directly from the definition \eqref{R}. Hence we can apply
 \eqref{relation-semi-norms} for $q:=-1$, and
 \eqref{semi-norms-sphere-estimate} for $g:=f$ to obtain
 \begin{eqnarray*}
 p_k(\cR[f]) & \overset{\eqref{relation-semi-norms}}{\le} & 
 (k+1)^{1+k}\varrho_k(\cR[f])\\
 & \overset{\eqref{semi-norms-sphere-estimate}}{\le} & (k+1)^{1+k}C(m,k)\varrho_k(f)\\
 & \overset{\eqref{relation-semi-norms}}{\le} & 
 (k+1)^{1+k}C(m,k)p_k(f).
 \end{eqnarray*}
\qed
Further restrictions to even smaller Fr\'echet spaces in the domain of $\cR$
and in the target space are necessary to invert $\cR$, since one can easily
see that the kernel of $\cR$ contains all odd functions. 
\begin{theorem}[$\cR$ as invertible mapping]\label{thm:inverse-radon}
The extended Radon transform $\cR$ restricted to the Fr\'echet
space of
even and positively $(-m)$-homogeneous functions of class $C^\infty(\R^{m+1}\setminus
\{0\})$ is a continuous and bijective linear mapping onto the Fr\'echet space of
even and positively $(-1)$-homogeneous smooth functions on $\R^{m+1}\setminus
\{0\})$ with a continuous inverse $\cT$.
\end{theorem}
\proof
The spherical Radon transform $\widehat{\cR}$ is injective on the space
of even continuous functions on the sphere $\S^m$; 
see \cite[Proposition 3.4.12]{groemer}. If $g_1,g_2$ are both, say $q$-homogeneous,
even, 
and smooth on $\R^{m+1}\setminus\{0\}$, then the identity
$$
\cR[g_1](Z)=\cR[g_2](Z)\quad\Foa Z\in\R^{m+1}\setminus\{0\}
$$
implies by definition (see \eqref{R})
$$
\widehat{\cR}[g_1|_{\S^m}](Z/|Z|)=
\widehat{\cR}[g_2|_{\S^m}](Z/|Z|)\quad\Foa Z\not= 0,
$$
from which we infer $g_1|_{S^m}=
g_2|_{S^m}$ since the  spherical Radon transform
$\widehat{\cR}$ is injective on such even functions. But then, by positive $q$-homogeneity,
$$
g_1(Z)=|Z|^{q}g_1|_{S^m}\Big(\frac{Z}{|Z|}\Big)=
|Z|^{q}g_2|_{S^m}\Big(\frac{Z}{|Z|}\Big)=g_2(Z).
$$
In addition, Helgason presents in \cite[Chapter III, Theorem 1.11]{helgason_radon}
an explicit formula for the preimage $\widehat{\cT}(\varphi)$ of any even
 function $\varphi$ of class $C^\infty(\S^m)$, which means
 that the spherical Radon transform is bijective on the space of
 even smooth functions on the sphere $\S^m$. Without using 
 Helgason's explicit expression
 for this inverse we can use its existence
 to define for a positively $(-1)$-homogeneous
 and even function $\Phi\in C^\infty(\R^{m+1}\setminus\{0\})$ the inverse
 $\cT:=\cR^{-1}$ of the extended Radon transformation $\cR$ by
 \begin{equation}\label{inverse-radon}
 \cT[\Phi](y):=\frac{1}{|y|^m}\widehat{\cT}[\Phi|_{\S^m}]\Big(\frac{y}{|y|}\Big).
 \end{equation}
Indeed, one calculates for an even positively $(-m)$-homogeneous
function $f\in C^\infty(\R^{m+1}\setminus\{0\})$, setting $\pi_{\S^m}(y):=y/|y|$,
\begin{eqnarray*}
\cT\big[\cR[f]\big](y) & = & \frac{1}{|y|^m}\hT\big[\big(\cR[f]\big)|_{\S^m}\big]\circ
\pi_{\S^m}(y)=\frac{1}{|y|^m}\hT\big[\Big((|\cdot|^{-1}\hR[f|_{\S^m}])\circ\pi_{\S^m}
\Big)|_{\S^m}\big]\circ
\pi_{\S^m}(y)\\
& = & \frac{1}{|y|^m}\hT\big[\Big(\hR[f|_{\S^m}]\circ\pi_{\S^m}\Big)|_{\S^m}
\big]\circ
\pi_{\S^m}(y)=\frac{1}{|y|^m}\hT\big[\Big(\hR[f|_{\S^m}]\Big)|_{\S^m}
\big]\circ
\pi_{\S^m}(y)\\
& = & \frac{1}{|y|^m}\hT\big[\hR[f|_{\S^m}]
\big]\circ
\pi_{\S^m}(y)=\frac{1}{|y|^m}f|_{\S^m}\circ\pi_{\S^m}(y)
 = 
\frac{1}{|y|^m}f(y/|y|) = f(y),
\end{eqnarray*}
where only in the very last equation we have used the positive
$(-m)$-homogeneity of $f.$ A similar calculation also shows
that
$\cR\big[\cT[\Phi]\big]=\Phi$ for all even, smooth, positively
$(-1)$-homogeneous functions $\Phi$ on $\R^{m+1}\setminus\{0\}.$
Since for any $q\in\R$ the space of even positively $q$-homogeneous
functions of class $C^\infty(\R^{m+1}\setminus\{0\})$ forms
a closed linear subspace of the Fr\'echet space $C^\infty(\R^{m+1}\setminus
\{0\})$ with respect to the family of seminorms $p_k$, $k=0,1,2,\ldots,$
we have established that $\cR$ is a continuous (see 
Corollary \ref{cor:radon-continuous}), bijective linear mapping 
from the Fr\'echet space of even positively $(-m)$-homogeneous smooth
mappings to the Fr\'echet space of even positively
$(-1)$-homogeneous smooth mappings on $\R^{m+1}\setminus\{0\}$. Hence we can
apply the open mapping theorem, in particular \cite[2.12 Corollaries (a), (b)]{rudin-fa-book_1973}, to obtain the continuity of $\cT$.
\qed

Now we have collected all properties needed to prove the parametric ellipticity
of the Cartan area integrand $\cA^F$ generated by a Finsler metric $F=F(x,y)$.
Let us recall 
this notion of convexity (see, e.g.,  \cite[p. 298]{HilvdM-dominance})
that is optimal for
any Cartan integrand $\cC=\cC(x,Z)$ satisfying the homogeneity condition
\begin{equation}\label{cartan-homogeneous}
\cC(x,tZ)=t\cC(x,Z)\quad\Foa (x,Z)\in \R^{m+1}\times (\R^{m+1}\setminus\{0\}),
\, t>0.
\end{equation}
by virtue of Euler's relation $\cC_{ZZ}(x,Z)Z=0.$
\begin{definition}[(Parametric) ellipticity]\label{def:ellipticity}
A Cartan integrand
$\cC=\cC(x,Z)\in C^2(\R^{m+1}\times\R^{m+1}\setminus\{0\})$ 
 satisfying \eqref{cartan-homogeneous}
is called
{\em
(parametric) elliptic} if and only if for every $R_0>0$ there is some number
$\lambda_{\cC}(R_0)>0$ such that the Hessian
$\cC_{ZZ}(x,Z)-\lambda_{\cC}(R_0)\cA^E_{ZZ}(Z)$
is positive semi-definite
for all $(x,Z)\in\overline{B_{R_0}(0)}
\times (\R^{m+1}\setminus\{0\}).$
\end{definition}
Recall from the remarks following Theorem \ref{thm:finsler-area}
that we denoted the Euclidean
metric by $E(y)=|y|$, which explains the notation $\cA^E$ for the Cartan area
integrand generated by $E$, and we have noticed there that $\cA^E(Z)=|Z|$.

\begin{theorem}[Ellipticity]\label{thm:ellipticity}
Let $F=F(x,y)$ be a reversible Finsler metric on $\R^{m+1}$. Then for 
each $x\in\R^{m+1}$
there exist constants $0<\lambda_1^F(x)\le\lambda_2^F(x)$ such that
\begin{equation}\label{ellipticity}
\lambda_1^F(x)\xi\cdot\cA^E_{ZZ}(Z)\xi\le
\xi\cdot\cA^F_{ZZ}(x,Z)\xi\le
\lambda_2^F(x)\xi\cdot\cA^E_{ZZ}(Z)\xi\quad\Foa \xi,Z\in\R^{m+1},\,Z\not=0.
\end{equation}
\end{theorem}
The Hessian $\cA^F_{ZZ}$ is $(-1)$-homogeneous in the $Z$-variable, and
for the Euclidean metric one computes
$$
|Z|\xi\cdot\cA^E_{ZZ}(Z)\xi=|\xi|^2-|Z|^{-2}(Z\cdot\xi)^2=|\pi_{Z^\perp}(\xi)|^2,
$$
where $\pi_{Z^\perp}$ denotes the orthogonal projection onto the subspace
$Z^\perp:=\{\eta\in\R^{m+1}:Z\cdot\eta=0\}.$ Therefore we deduce from
the left inequality  in \eqref{ellipticity}  by homogeneity
$$
\inf_{\bar{Z}\in\S^m,\atop\bar{\eta}\in\S^m\cap \bar{Z}^\perp}
\bar{\eta}\cdot\cA^F_{ZZ}(x,\bar{Z})\bar{\eta}=\inf_{Z\not=0,\atop \eta\in Z^\perp\setminus\{0\}}\frac{|Z|\eta\cdot\cA^F_{ZZ}(x,Z)\eta
}{|\eta|^2}\ge \lambda_1^F(x)>0.
$$
Since the right-hand side is positive for each $x\in\R^{m+1}$ according to
Theorem \ref{thm:ellipticity}, and the
left-hand side is (Lipschitz) continuous in $x$, 
it attains its positive minimum $\lambda(R_0)>0$
on $\overline{B_{R_0}(0)}\subset\R^{m+1}$. This implies
$$
|Z|\eta\cdot\cA^F_{ZZ}(x,Z)\eta\ge\lambda(R_0)|\eta|^2\quad\Foa Z\not=0,\,x\in\overline{B_{R_0}(0)},\,\eta\in Z^\perp,
$$
which readily translates to
\begin{equation}\label{parametric-ellipticity}
|Z|\xi\cdot\cA^F_{ZZ}(x,Z)\xi\ge\lambda(R_0)|\pi_{Z^\perp}(\xi)|^2\quad\Foa Z\not=0,\,
x\in\overline{B_{R_0}(0)},\,\xi\in\R^{m+1}.
\end{equation}
Combining this with Lemma \ref{lem:equal-Cartan}
we have shown that Theorem \ref{thm:ellipticity} implies the following:
\begin{corollary}[Ellipticity]\label{cor:ellipticity}
Let $F=F(x,y)$ be a Finsler metric on $\R^{m+1}$ satisfying assumption
{\rm (GA)}. Then the corresponding Cartan area integrand $\cA^F$ is
(parametric) elliptic in the sense of Definition \ref{def:ellipticity}.
\end{corollary}

{\sc Proof of Theorem \ref{thm:ellipticity}:}
Fix $x\in\R^{m+1}$. We claim that the function
$$
Z\mapsto \cA^F(x,Z)-\frac 1n |Z|
$$ 
is convex for each $n\ge n_0$ if $n_0=n_0(x)$ is sufficiently large. From this
statement the theorem follows.

To prove the claim we first apply Lemma \ref{lem:pointwise-bounds}
to the Finsler structures $F_1(x,y):=E(y)=|y|$ and $F_2(x,y):=F(x,y)$
to deduce from
$$
c_1(x)E(y):=\min_{\eta\in\S^m}F(x,\eta)|y|\le F(x,y)\le \max_{\eta\in\S^m}F(x,\eta)|y|=:
c_2(x)E(y)
$$
the corresponding bounds on the Cartan area integrands $\cA^E$ and $\cA^F$,
\begin{equation}\label{A-bounds2}
c_1^m(x)\cA^E(Z)\le\cA^F(x,Z)\le c_2^m(x)\cA^E(Z)\quad\Foa Z\in\R^{m+1}.
\end{equation}
Notice that $0<c_1(x)\le c_2(x)$ since $F$ as a Finsler metric with the properties 
(F1)
and (F2) satisfies $F(x,y)>0$ for $y\not= 0$; 
see, e.g., \cite[Theorem 1.2.2]{bao-chern-shen_2000}.
Thus, for all $n\in\N$ with $n\ge  n_1(x)$, where $n_1(x)$ is the smallest integer greater or equal $c_1(x)^{-1},$ the expression
$$
\Phi_n(x,Z):=\cA^F(x,Z)-\frac 1n |Z|
$$
is an even,
\heikodetail{

\bigskip

\tt\xx CHECK WHY ?? NEEDS TO BE MENTIONED EARLIER, PROBABLY
VOR LEMMA \ref{lem:pointwise-bounds} ODER ETWAS SPAETER, Z.B.
DIREKT NACH LEMMA \ref{lem:cartan-radon}
$F$ is reversible\xx\xx

\yy Indeed $F$ has to be reversible at this stage for the proof to work, as it has to lie in the image of the inverse Radon transform in some sense. 

\bigskip

} 
positively $1$-homogenous smooth function on $\R^{m+1}\setminus\{0\}$
(recall that $x$ is fixed), with $\Phi_n(x,Z)>0$ for $Z\not= 0;$ 
cf. Lemma \ref{lem:area-comparison} and our 
remarks directly following Theorem \ref{thm:finsler-area}.
Hence, by virtue of Theorem \ref{thm:inverse-radon}, the inverse $\cT$
of the extended Radon transfrom $\cR$ can be applied to $1/\Phi_n(x,\cdot)$.

Before doing so
we 
use Lemma \ref{lem:cartan-radon} to rewrite the function $\Phi_n$
as
$$
\Phi_n(x,Z)=\cA^F(x,Z)-\frac 1n |Z|  =  \frac{1}{\cR[F(x,\cdot)^{-m}]}- \frac 1n 
\frac{1}{\cR[|\cdot|^{-m}]}=\frac{\cR[|\cdot|^{-m}]-\frac 1n \cR[F(x,\cdot)^{-m}]}{
\cR[F(x,\cdot)^{-m}]\cR[|\cdot|^{-m}]}.
$$
Applying $\cT$ to the function
 $1/\Phi_n(x,\cdot)$ (for fixed $x$ and $n\ge n_1(x)$) now yields
\begin{equation}\label{inv-radon-Phi_n} 
 \cT\Big(\frac{1}{\Phi_n(x,\cdot)}\Big)=
\cT\Big(\frac{\cR[F^{-m}(x,\cdot)]\cR[|\cdot|^{-m}]}{
\cR[|\cdot|^{-m}]-\frac 1n \cR[F^{-m}(x,\cdot)]}\Big).
\end{equation} 
The argument of $\cT$ on the right-hand side is an even $(-1)$-homogenous 
function that tends to the expression $\cR[F^{-m}(x,\cdot)]$ as $n\to\infty$
in the Fr\'echet space of even $(-1)$-homogeneous smooth functions on
$\R^{m+1}\setminus\{0\}$ with respect to the 
topology given by the family of seminorms $p_k$, $k=0,1,2,\ldots,$
introduced in \eqref{semi-norms-space}.
By the continuity of $\cT$ with respect to that convergence (granted
by Theorem \ref{thm:inverse-radon}) we find that
\begin{equation}\label{limit}
\Psi_n(x,\cdot):=\Big[\cT\Big(\frac{1}{\Phi_n(x,\cdot)}\Big)\Big]^{-1/m}\longrightarrow
F(x,\cdot)\quad\textnormal{as $n\to\infty$,}
\end{equation}
in the Fr\'echet space of even $1$-homogeneous smooth functions on $\R^{m+1}\setminus
\{0\}$ with respect to the family of seminorms $p_k$, $k=0,1,2\ldots$ .
(Notice that the left-hand side of \eqref{inv-radon-Phi_n} is even and
positively $(-m)$-homogeneous
according to Theorem \ref{thm:inverse-radon}, so that
the left-hand side of \eqref{limit} is even and positively $1$-homogeneous.)

$F$ itself is a reversible Finsler metric satisfying (F1) and (F2) by assumption.
In particular, for the positively $0$-homogeneous Hessian
$(F^2/2)_{yy}(x,y)$ one finds for fixed $x$ a constant $\Lambda_1^F(x)>0$ such that
\begin{equation}\label{f2}
\xi\cdot \Big(\frac{F^2}{2}\Big)_{yy}(x,y)\,\xi\ge \Lambda_1^F(x)|\xi|^2\quad\Foa \xi\in\R^{m+1},\,y\not=0.
\end{equation}
\heikodetail{

\bigskip

Since $(F^2/2)_{yy}$ is positively $0$-homogeneous and by (F2) one finds
\begin{equation}\label{detail-convex}
\frac{\xi}{|\xi|}\cdot \Big(\frac{F^2}{2}\Big)_{yy}(x,\frac{y}{|y|})\frac{\xi}{|\xi|^2}=
\frac{\xi\cdot (F^2/2)_{yy}(x,y/|y|)\xi}{|\xi|^2}>0,
\end{equation}
whence
$$
\inf_{\zeta,\eta\in\S^m}\zeta\cdot (F^2/2)_{yy}(x,\eta)\zeta=\Lambda_1^F(x)>0.
$$
Otherwise there existed sequences $\zeta_i\to\bar{\zeta}$ and $\eta_i\to
\bar{\eta}$ in $\S^m$ such that
$$
0\le\zeta_i\cdot (F^2/2)_{yy}(x,\eta_i)\zeta_i<\frac 1i,
$$
which in the limit gives $\bar{\zeta}\cdot (F^2/2)_{yy}(x,\bar{\eta})\bar{\zeta}=0$
contradicting \eqref{detail-convex}.

\bigskip

}
By virtue of the convergence in \eqref{limit} there is some $n_0(x)\ge n_1(x)$
such that
\begin{equation}\label{f2-limit}
\xi\cdot\Psi_n(x,y)\,\xi\ge
\frac{\Lambda_1^F}{2}|\xi|^2\quad\Foa n\ge n_0(x),\,\xi\in\R^{m+1},\,y\not=0,
\end{equation}
which qualifies 
$$
\Psi_n(x,\cdot)=\Big[\cT\Big(\frac{1}{\Phi_n(x,\cdot)}\Big)\Big]^{-1/m}
$$
 as a reversible Finsler metric on $\R^{m+1}$ for each $n\ge n_0(x)$. Consequently, by Busemann's convexity result, 
Theorem \ref{thm:busemann}, the corresponding Cartan area integrand
$
\cA^{\Psi_n(x,\cdot)}=\cA^{\Psi_n(x,\cdot)}(x,Z)$ is convex in the $Z$-variable.
By means of Lemma \ref{lem:cartan-radon} we can rewrite this Cartan area integrand
as
$$
\cA^{\Psi_n(x,\cdot)}(x,Z)=\frac{1}{\cR[\Psi_n^{-m}(x,\cdot)](Z)}=
\frac{1}{\cR [\cT\Big(\frac{1}{\Phi_n(x,\cdot)}\Big)]}=\Phi_n(x,Z)=\cA^F(x,Z)-\frac
1n |Z|,
$$
which establishes the left inequality in \eqref{ellipticity} for $\lambda_1^F(x)=1/n_0(x).$
The right inequality in \eqref{ellipticity} simply follows from the fact that
$\cA^F(x,Z)$ is smooth on $\R^{m+1}\setminus\{0\}$ and positively $1$-homogeneous.
\heikodetail{

\bigskip

Indeed, $\xi=\alpha\eta+\beta Z$  with $|\eta|=1$ and $\eta \cdot Z= 0$  so that
\begin{eqnarray*}
\xi\cdot |Z|\cA^F_{ZZ}(x,Z)\xi & = & \alpha^2\eta 
\cdot |Z|\cA^F_{'ZZ}(x,Z)\eta\le\alpha^2 C(x)\\
& = & C(x)\alpha\eta\cdot  |Z|\cA^E_{ZZ}(x,Z)(\alpha\eta)=\xi\cdot
|Z|\cA^E_{ZZ}(x,Z)\xi
\end{eqnarray*}
for 
$$
C(x):=\sup_{\zeta,\mu\in\S^m}\mu\cdot \cA^F_{ZZ}(x,\zeta)\mu.
$$

\bigskip

}

\qed

\setnumbers
\section{Proofs of the main results}\label{sec:3}
{\bf Finsler-minimal graphs.}\, First, we will combine our results of Section \ref{sec:2}
with well-known results of Jenkins \cite{jenkins_1961b}, L. Simon \cite{simon_1977a},
and Winklmann \cite{winklmann-bernstein} on solutions of the non-parametric
Euler-Lagrange equations of elliptic Cartan functionals (whose integrand does
not depend explicitly on the position vector), to prove the Bernstein result.

\medskip

{\sc Proof of Theorem \ref{thm:bernstein}:}\,
Let $F=F(y)$ be a Minkowski metric on $\R^{m+1}$.
Any entire Finsler-minimal graph $\{(u,f(u)):u\in\R^m\}$ is by definition
a critical immersion of Finsler area, which -- according to
Theorem \ref{thm:finsler-area} -- can be written as the variational integral
\begin{equation}\label{cartan-functional}
\int_{\R^m}\cA^F((X_{u^1}\wedge\ldots\wedge X_{u^m})(u))\,du^1\wedge\ldots
\wedge du^m,
\end{equation}
where $X(u)=(u,f(u))$ for $u\in\R^m$. Here, the Cartan area
integrand $\cA^F=\cA^F(Z)$ does not depend on the position vector,  
is smooth on $\R^{m+1}\setminus\{0\}$ and positively
$1$-homogeneous (see \eqref{H} and our remarks directly following Theorem \ref{thm:finsler-area}). In addition, since $F=F(y)$ does
not depend on the $x$-variable, we have the simple estimate (cf. \eqref{definite})
\begin{equation}\label{F-minkowski-bounds}
m_F|y|=\min_{\eta\in\S^m}F(\eta)|y|\le F(y/|y|)|y|\le\max_{\zeta\in\S^m}F(\zeta)|y|=M_F|y|
\end{equation}
with $0<m_F\le M_F$ by the defining properties (F1) and (F2) which together imply
that $F(y)>0$ as long as $y\not=0.$ Lemma \ref{lem:pointwise-bounds} implies
that 
\begin{equation}\label{AF-minkowski-bounds}
m_F^m|Z|\le \cA^F(Z)\le M_F^m|Z|.
\end{equation}
By virtue of the general assumption (GA) and Corollary \ref{cor:ellipticity}
we know that the Cartan area integrand $\cA^F=\cA^F(Z)$ is (parametric) elliptic.
Since there is no $x$-dependence here, we therefore have the simplified
ellipticity condition (cf. \eqref{parametric-ellipticity})
\begin{equation}\label{AF-minkowski-ellipticity}
|Z|\xi\cdot\cA^F_{ZZ}(Z)\xi\ge\lambda|\pi_{Z^\perp}(\xi)|^2
\end{equation}
for some positive constant $\lambda$.

By scaling we may assume that $\mu:=\min\{m_F,\lambda\}\ge 1$, 
since a critical  immersion $X$ for
\eqref{cartan-functional} is also critical for this integral with $\cA^F$
replaced by $\mu^{-1}\cA^F$.

Thus we have verified that the general assumptions (i)--(iii) on the
Cartan integrand in \cite[p. 266]{simon_1977a} (or alternatively without
the need to scale in \cite[p. 181]{jenkins_1961b}) are satisfied. The non-parametric
Euler-Lagrange equation is given, e.g., in \cite[formula (1)]{simon_1977a},
and the function $f\in C^2(\R^m)$ generating the entire
Finsler-minimal graph $\{(u,f(u)):u\in\R^m\}$ solves that equation on $\R^m$.

The Bernstein theorem for $m=2$ formulated in part (i) now follows either from 
\cite[Theorem 3 \& Corollary]{jenkins_1961b} or from the first
statement in \cite[Corollary 1]{simon_1977a}. (Notice for the latter that
Simon introduces the class $\mathscr{M}'$ of generalized surfaces that 
can be represented as a graph of a $C^2$-function, and this class is contained
in the class $\mathscr{M}$ of hypersurfaces 
mentioned in the first statement of his Corollary 1.)
For $m=3$ we refer to the explicit statement regarding entire solutions of the
non-parametric Euler-Lagrange equation in the second part of 
\cite[Corollary 1]{simon_1977a}.

For the proof of part (ii) of Theorem \ref{thm:bernstein} we recall the seminorms
$\varrho_k$ on the sphere introduced in \eqref{semi-norms-sphere}, and assume
to the contrary that for 
a sequence of Minkowski norms $F_n=F_n(y)$ on $\R^{m+1}$ converging up to third
order to the Euclidean metric $E(y)=|y|$,
\begin{equation}\label{convergence-third-order}
\varrho_3(F_n-E)\longrightarrow 0\quad\textnormal{as $n\to\infty,$}
\end{equation}
we have entire
Finsler-minimal graphs $X_n(u):=(u,f_n(u))$ that are not affine planes. Here,
$X_n$ is critical for the Finsler area $\area_{\R^m}^{F_n}.$
From \eqref{convergence-third-order} we deduce, in particular,  by means of the triangle inequality, that 
there is $n_0\in\N$ such that
$$
\frac 12 =\frac 12 \varrho_0(E)\le\varrho_0(F_n)\le 2 \varrho_0(E)=2\quad\Foa n\ge n_0.
$$
Hence, by the chain rule and by \eqref{convergence-third-order},
\begin{equation}\label{convergence2}
\varrho_3(F_n^{-m}-E^{-m})\longrightarrow 0\quad\textnormal{as $n\to\infty,$}
\end{equation}
which according to Lemma \ref{lem:semi-norms-sphere} implies for the (linear)
extended Radon transformation $\cR$:
\begin{equation}\label{Rconvergence2}
\varrho_3(\cR[F_n^{-m}]-\cR[E^{-m}])=
\varrho_3(\cR[F_n^{-m}-E^{-m}])\longrightarrow 0\quad\textnormal{as $n\to\infty.$}
\end{equation}
Again, by the triangle inequality,
\begin{equation}\label{Rconv-estimate}
\frac 12 \le C_1(m):=\frac 12 \varrho_3(\cR[E^{-m}])\le \varrho_3(\cR[F_n^{-m}])\le 2\varrho_3(
\cR[E^{-m}])=:C_2(m)
\quad\textnormal{for all $n\ge n_1$}
\end{equation}
for some $n_0\le n_1\in\N.$ Notice that 
$$
2C_1(m)\ge\varrho_0(\cR[E^{-m}])=
\max\{(\cA^E(Z))^{-1}=|Z|^{-1},\, Z\in\S^m\}=1,
$$
where we used Lemma \ref{lem:cartan-radon}.
By the same lemma, on the other hand,
for $Z\in\R^{m+1}\setminus\{0\}$, 
\begin{equation}\label{Rterm}
\cA^{F_n}(Z)-\cA^E(Z)  \overset{\eqref{cartan-radon}}{=}  \frac{1}{\cR[F_n^{-m}](Z)}-
\frac{1}{\cR[|\cdot|^{-m}](Z)}=\frac{\cR[|\cdot|^{-m}](Z)-\cR[F_n^{-m}](Z)}{
\cR[F_n^{-m}](Z)\cR[|\cdot|^{-m}](Z)}.
\end{equation}
Taking the $\varrho_3$-seminorm of this expression with a careful application
of the  Leibniz product rule, one can then use \eqref{Rconv-estimate} to find
a constant $C_3(m)\ge 0 $ such that 
$$
\varrho_3(\cA^{F_n}(Z)-\cA^E(Z))\le C_3(m)\varrho_3(\cR[|\cdot|^{-m}]-\cR[F_n^{-m}])
\quad\textnormal{for all $n\ge n_1$},
$$
and the right-hand side converges to zero as $n\to \infty$ by virtue of \eqref{Rconvergence2}.
Thus, for any given $\eta >0$ there is $n_2\in\N$ with $n_2\ge n_1$, such that
$$
\sum_{|\alpha|\le 3}|D^\alpha (\cA^{F_n}-\cA^E)(\xi)|<\eta\quad\Foa\xi\in\S^m,\,n\ge n_2,
$$
which is exactly the condition that 
Simon requires for his Bernstein result for solutions
of the non-parametric Euler-Lagrange equations of \eqref{cartan-functional}
in dimensions $m\le 7$ (\cite[Corollary 2 \& (6)]{simon_1977a}.
Hence, all Finsler-minimal graphs $X_n$, for $n\ge n_2$, are 
affine planes, contradicting our assumption.

It remains to prove part (iii) of Theorem \ref{thm:bernstein}. For this 
we verify the conditions Winklmann assumes for his Bernstein result
\cite[Theorem 4.1]{winklmann-bernstein} for  $C^2$-solutions of the 
non-parametric Euler-Lagrange equation. Indeed, the present Cartan area integrand
$\cA^F=\cA^F(Z)$ satisfies the positive $1$-homogeneity and ellipticity condition
required in \cite[(2.1) \& (2.2)]{winklmann-bernstein} (see \eqref{H} and
\eqref{AF-minkowski-ellipticity} in our context),
and the $C^3$-vicinity
quantized with a constant $\delta_\star(m,\gamma)$ in \cite[p. 383]{winklmann-bernstein}
translates to condition \eqref{C3-close} by means of the indirect argument 
presented for the proof of part (ii), where the final contradiction will be
obtained by \cite[Theorem 4.1]{winklmann-bernstein}  in this case,
which concludes the proof.
\qed

\bigskip

The uniqueness result for Finsler-minimal graphs, Theorem \ref{thm:uniqueness},
will be established connecting our results of Section \ref{sec:2} to a more
recent weighted energy estimate of Hildebrandt and Sauvigny 
\cite[Theorem 3.1]{Hil-sauvigny-energy-estimate-2009}.

{\sc Proof of Theorem \ref{thm:uniqueness}:}\,
In the proof of Theorem \ref{thm:bernstein} we have established
already that both Finsler-minimal graphs $f_1,f_2\in C^0(\bar{\Omega}\setminus K)
\cap C^2(\Omega\setminus K)$ are solutions of the non-parametric Euler-Lagrange
equation of the Cartan functional \eqref{cartan-functional} on the domain
$\Omega\setminus K\subset\R^m$. The ellipticity condition 
\eqref{AF-minkowski-ellipticity} is a simpler version of the ellipticity condition
\cite[(3.7)]{Hil-sauvigny-energy-estimate-2009} since we have no explicit
dependence on the position vector in the Cartan area integrand in 
\eqref{cartan-functional}. This translates into an ellipticity condition
of the non-parametric functional
(where we can ignore the $x$- and $z$-dependence in 
\cite[(3.14)]{Hil-sauvigny-energy-estimate-2009}),
which is explicitly assumed in \cite[Theorem 3.1]{Hil-sauvigny-energy-estimate-2009}.
Hence we can apply this theorem in the simpler situation of the homogeneous
non-parametric Euler-Lagrange equation, i.e., for the right-hand side $H=H(x,z)\equiv 0$
\heikodetail{

\bigskip

Herr Overath erklaert:
\yy Die pde \cite[(4.16)]{Hil-sauvigny-energy-estimate-2009} fuer einen kritischen Graphen einer Funktion $f$ schreibt sich hier fuer den nichtparametrische Integrand $I(\tilde{x},z,p):= 
A^F(\tilde{x},z, -p, 1)$ als
$$
\mathrm{div}I_p(\cdot,f,\nabla f) = m H(\cdot, f).
$$
Darin ist dann $H(\tilde{x},z,p)=I_z(\tilde{x},z,p)
=A^F_{x^3}(x,Z)$ mit $x=(x^1,x^2,x^3)=(\tilde{x},z), 
\tilde{x}=(x^1,x^2),Z=(-p^1,-p^2,1)$. Da wir $F$ als Minkowski-Metrik 
gewaehlt haben, folgt $H\equiv 0$! Es wuerde aber auch schon genuegen $F$ als unabhaengig von der 3ten Komponente $x^3$ von $x$ zu waehlen.

\bigskip

\xx\xx Ist das eine INTERESSANTE Klasse von Finsler-Metriken?

\bigskip

\xx PLEASE EXPLAIN - seems to be valuable remark! \xx \yy I just intended to explain how one computes the quantity $H$ of \cite{Hil-sauvigny-energy-estimate-2009} in our setting. \yy

\bigskip

} 
in \cite[(3.22)]{Hil-sauvigny-energy-estimate-2009}. The resulting weighted
energy estimate of Hildebrandt and Sauvigny reads therefore in our situation as follows:
\begin{equation}\label{energy-estimate}
\int_{\Omega\setminus K}\mu(f_1,f_2)|\nabla f_1-\nabla f_2|^2\,dx\le
\frac{2}{\lambda}\int_{\partial\Omega}|f_1-f_2|\,d\mathscr{H}^{m-1},
\end{equation}
where $\mu(f_1,f_2)(u):=\big(\max\{\sqrt{1+\nabla f_1^2(u)},
\sqrt{1+\nabla f_2^2(u)}\}\big)^{-3}$ for $u\in\Omega\setminus K.$
Now we can proceed as in the proof of 
\cite[Theorem 4.1]{Hil-sauvigny-energy-estimate-2009}. Assuming that $f_1|_{\partial\Omega}=
f_2|_{\partial\Omega}$ we conclude $\nabla f_1\equiv \nabla f_2$ on $\Omega
\setminus K$, and since $\Omega\setminus K$ is a domain this implies
that $f_1-f_2\equiv const$ on $\Omega\setminus K$. But $f_1-f_2$ is of class
$C^0(\bar{\Omega}\setminus K)$ and $K$ is compactly contained in $\Omega$, so that
$f_1\equiv f_2$ on $\partial \Omega$ leads to equality of $f_1$ and $f_2$ on all 
of $\bar{\Omega}\setminus K.$
\qed

\bigskip

Uniqueness results like Theorem \ref{thm:uniqueness} can be used to
extend solutions $f_1$ of certain partial differential equations on domains
$\Omega\setminus K$ to all of $\Omega$, if the Dirichlet problem on $\Omega$
with prescribed boundary data $f_1$ on $\partial\Omega$ can be solved by a function,
say $f_2$. The uniqueness theorem then implies that this solution $f_2$ coincides
with $f_1$ on $\Omega\setminus K$, so that
the original singular set $K$ of $f_1$ is
removed by extending $f_1$ to the solution $f_2$ on all of $\Omega.$
Such a result can be found in \cite[Remark 4.7]{Hil-sauvigny-energy-estimate-2009}
but under slightly stronger assumptions on the singular set $K$, so we will draw
from L. Simon's contribution \cite{simon_1977b} to prove Theorem \ref{thm:removability}.

{\sc Proof of Theorem \ref{thm:removability}.}\,
As in the proof of Theorem \ref{thm:bernstein} we find that any Finsler-minimal
graph $\{(u,f(u)):u\in \Omega\setminus K\}$ of class $C^2$ 
in Finsler-Minkowski space $(\R^{m+1},F=F(y))$ 
is generated by a function $f\in C^2(\Omega\setminus K)$
that solves the non-parametric Euler-Lagrange equation for the variational
integral \eqref{cartan-functional} as, e.g., 
presented in \cite[formula (5)]{simon_1977b}. Also the general assumptions
\cite[formulae (1)--(4)]{simon_1977b} are satisfied in our present situation
(cf. \eqref{AF-minkowski-bounds} and \eqref{AF-minkowski-ellipticity}, modulo
the simple scaling argument mentioned in the proof of Theorem \ref{thm:bernstein}
to obtain $\mu:=\min\{m_F,\lambda\}\ge 1$ for the exact  factors
in inequalities
\cite[(2) \& (3)]{simon_1977b}). Now, Theorem \ref{thm:removability} follows
from \cite[Theorem 1]{simon_1977b}.
\qed

\bigskip

{\bf Global results for Finsler-minimal immersions.}\,
Combining our results of Section \ref{sec:2} with global
results on critical immersions of Cartan functionals by Clarenz et. al.,and
Winklmann 
\cite{clarenz_enclosure,clarenz-vdM-isop,winklmann-isop,winklmann-existence-uniqueness} we establish in the following the convex hull
property formulated in Theorem \ref{thm:convex-hull}, the isoperimetric
inequalities in Theorem \ref{thm:isop}, and the existence results for
Finsler-minimal immersions in Theorem \ref{thm:existence}.

{\sc Proof of Theorem \ref{thm:convex-hull}.}\,
If $F=F(y)$ is a  Minkowski metric on $\R^{m+1}$ then $F^{-m}$
\heikodetail{

\bigskip

\tt \yy $(\R^{m+1},F)$ is a Minkowski space

\bigskip

PLEASE EXPLAIN\xx

\yy Siehe meine Bemerkung bei der Definition zur Minkowski metric eingangs der Arbeit. \yy
} 
is smooth
on $\R^{m+1}\setminus \{0\}$
by virtue of \eqref{F-minkowski-bounds}, so that Lemma \ref{lem:cartan-radon}
in combination with Corollary \ref{cor:radon-continuous} leads to a
Cartan area integrand $\cA^F=\cA^F(Z)$ that is positively $1$-homogeneous
and smooth on $\R^{m+1}\setminus\{0\}.$ In addition, we use the general assumption
(GA) together with Corollary \ref{cor:ellipticity} as in the proof of Theorem 
\ref{thm:bernstein} to show that $\cA^F$ is (parametric) elliptic with the uniform
estimate \eqref{AF-minkowski-ellipticity}, which coincides with 
the ellipticity condition
used in \cite[Definition 2.2]{clarenz_enclosure}. 
Since we assume that $X\in C^2(\mathscr{M},\R^{m+1})\cap C^0(\bar{\mathscr{M}},\R^{m+1})$
is a Finsler-minimal immersion we know that $X$ is critical for the Cartan
area functional \eqref{cartan-functional}; in the language of Clarenz \cite{clarenz_enclosure} $X$ is extremal for \eqref{cartan-functional}. This implies by
\cite[Theorem 1]{clarenz_enclosure} that its $\cA^F$-mean curvature $H_{\cA^F}$
vanishes since the Cartan area integrand $\cA^F$ does not depend on the position
vector. Hence all conditions of \cite[Theorem 2.3]{clarenz_enclosure}
are satisfied which establishes the convex hull property for $X$.
\qed

\bigskip

The first variation formula for general Cartan functionals derived by
R\"awer \cite{raewer-phd_1993}
\heikodetail{

\bigskip

Herr Overath schreibt:
\tt\yy aAuf diese Quelle habe ich momentan keinen Zugriff

BUT LET US MENTION RAEWER ANYWAY\xx

\bigskip

} 
and Clarenz \cite{clarenz-phd_1999} were used
by Clarenz et al. and Winklmann to derive various isoperimetric inequalities
for critical immersions of such functionals 
\cite{clarenz-vdM-compact_2001,clarenz-vdM-isop,winklmann-isop}. Some of these
results are used in the present context in connection with Section \ref{sec:2}
to derive corresponding results for Finsler-minimal immersions.

{\sc Proof of Theorem \ref{thm:isop}.}\,
The first two parts of this theorem deal with Finsler-minimal immersions in
Minkowski space. The Cartan area integrand $\cA^F$ derived from the Minkowski
metric $F=F(y)$ depends only on the $Z$-variable, and such Cartan functionals
without dependence on the position vector have been considered by Winklmann
in \cite{winklmann-isop}. Notice that Winklmann does not use positivity
of his Cartan integrand, although we have it here for $\cA^F(Z)$ automatically
by means of \eqref{AF-minkowski-bounds}. Now Winklmann's isoperimetric
inequality in \cite[Corollary 2.3]{winklmann-isop}
reads in our context as
\begin{equation}\label{winkl-isop}
\area_\mathscr{M}^F(X)\le \frac 1m \frac{R}{m_F} \|\cA^F_Z\|_{C^0(\S^m,\R^{m+1})}
\int_{\partial\mathscr{M}}dS.
\end{equation}
\heikodetail{

\bigskip

Ich hatte gefragt:
\tt \xx Was ist bei Winklmann gemeint mit $L:=\int_{\partial \mathscr{M}} dS$? 
Tatsaechlich das $(m-1)$-dimensuionale Hausdorffvolumen 
von $X(\partial \mathscr{M})$ oder der $(m-1)$-dimensionale 
Flaecheninhalt von $\partial \mathscr{M}$ unter der Immersion $X$? 
Vielleicht koennen wir dies nochmal diskutieren. 

Herr Overaths Antwort:
Mir erscheint $\area^{|\cdot|}_{\partial \mathscr{M}}(X|_{\partial \mathscr{M}})$ 
anstatt $\mathscr{H}^{m-1}(X(\partial\mathscr{M}))$ als passender, 
da dieser Ausdruck auch fuer nicht injektive $X$ mit $L$ 
uebereinstimmt. 

$dS$ wird im Winklmann-Kontext ueber die auf 
$\partial \mathscr{M}$ zurueckgeholte euklidische Metrik gebildet.

\bigskip

}
Recall that $R>0$ is the radius of the closed Finsler ball 
$\overline{B_R^F(0)}\subset\R^{m+1}$
containing the boundary $X(\partial\mathscr{M})$, and therefore by the convex
hull property, Theorem \ref{thm:convex-hull}, also $X(\bar{\mathscr{M}})$. 
Notice further, that $\overline{B_R^F(0)} \subset \overline{B_{R/m_F}(0)}$, which 
can be shown by an argument similar to the  one leading to
\eqref{proof:area-comparison-Finsler-unit-ball-inclusion} in the proof
of Lemma \ref{lem:area-comparison}.
To estimate the supremum-norm of $\cA^F_Z$ purely in terms of quantities
explicitly given by the Minkowski metric $F$ we use Lemma \ref{lem:cartan-radon}
and the differentiation rule
for the Radon transform in the simple form \eqref{diff-rule} (contracted by 
multiplication with $Z^\tau$ and summation over $\tau$ to compute in
a first step for $Z\in\S^m$  
\begin{eqnarray*}
\cA^F_{Z_\sigma}(Z)& \overset{\eqref{cartan-radon}}{=} &  \Big(\frac{1}{\cR[F^{-m}]}
\Big)_{Z_\sigma}(Z)=-\frac{(\cR[F^{-m}])_{Z_\sigma}(Z)}{\cR^2[F^{-m}](Z)}\\
 &\overset{\eqref{diff-rule}}{=}& \frac{1}{\cR^2[F^{-m}](Z)}\sum_{\tau=1}^{m+1}\cR\Big[Z^\tau\frac{\partial}{
\partial y^\tau}(y^\sigma F^{-m})\Big](Z)\\
& = & \frac{1}{\cR^2[F^{-m}](Z)}\sum_{\tau=1}^{m+1}\cR\Big[Z^\tau (\delta^\sigma_\tau
F^{-m}-my^\sigma F^{-m-1}F_{y^\tau})\Big](Z)\\
& = & \frac{Z^\sigma}{\cR[F^{-m}](Z)}-\frac{m}{\cR^2[F^{-m}](Z)}\sum_{\tau=1}^{m+1}
\cR[Z^\tau y^\sigma F^{-m-1}F_{y^\tau}](Z),
\end{eqnarray*}
\heikodetail{

\bigskip

\yy Remember that we can neglect the Euclidean norm $|Z|$ in the above computations since $|Z|=1$ as $Z\in\S^{m}$. \yy

\bigskip

PLEASE EXPLAIN\xx \yy I just was wondering why some factors involving $|Z|$ weren't present in the above computations until I saw that this comes from the choice of $Z$.\yy

}
where we used the linearity of the Radon several times. Therefore, we obtain
for the norm of the gradient $\cA^F_Z(Z)$  at $Z\in\S^m$ (again using linearity
of $\cR$)
\begin{eqnarray}
|\cA^F_Z(Z)|^2 & = & \frac{1}{\cR^2[F^{-m}](Z)}-\frac{2m}{\cR^3[F^{-m}](Z)}
\sum_{\sigma,\tau=1}^{m+1}\cR[Z^\tau Z^\sigma y^\sigma F^{-m-1}F_{y\tau}](Z)
\label{AF-estimate}\\
& & + \frac{m^2}{\cR^4[F^{-m}](Z)}\sum_{\sigma=1}^{m+1}\Big(\sum_{\tau=1}^{m+1}\cR[
Z^\tau y^\sigma F^{-m-1}F_{y^\tau}](Z)\sum_{\rho=1}^{m+1}\cR[Z^\rho y^\sigma
F^{-m-1}F_{y^\rho}](Z)\Big)\notag
\end{eqnarray}
Since the spherical Radon transform $\widehat{\cR}$ coincides with the extended
Radon transform $\cR$ on the sphere $\S^m$ we can use the explicit formula
\eqref{spherical-radon} to show that the second term on the right-hand side
above vanishes:
for $g(y):=\sum_{\tau=1}^{m+1}Z^\tau F^{-m-1}(y)F_{y^\tau}(y)$ we compute
for a fixed $Z\in\S^m$ by linearity of $\cR$
$$
\sum_{\sigma=1}^{m+1}\cR[Z^\sigma y^\sigma g(\cdot)](Z)=
\cR[(Z\cdot y) g(\cdot)]\overset{\eqref{spherical-radon}}{=}
\frac{1}{\mathscr{H}^{m-1}(\S^{m-1})}\int_{y\in\S^m\cap Z^\perp}(Z\cdot y) g(y)\,
d\H^{m-1}(y)=0.
$$
Consequently, it suffices to estimate the first and the last term on the right-hand
side of \eqref{AF-estimate} in terms of quantities related to the Minkowski
metric $F$.
By Lemma \ref{lem:cartan-radon} we immediately see that the first term coincides
with $(\cA^F(Z))^2$. For the last term in \eqref{AF-estimate} 
we set $h(y):=(Z\cdot F_y(y))/F(y)$ and use again the explicit
integral formula \eqref{spherical-radon} to estimate  for fixed $\sigma\in
\{1,\ldots,m+1\}$
\begin{eqnarray*}
\sum_{\tau=1}^{m+1}\cR[
Z^\tau y^\sigma F^{-m-1}F_{y^\tau}](Z) &= &\cR[y^\sigma (Z\cdot F_y(y))/F(y)F^{-m}(y)]\\
& = &
\frac{1}{\H^{m-1}(\S^{m-1})}\int_{y\in S^m\cap Z^\perp} y^\sigma 
h(y)F^{-m}(y)\,d\H^{m-1}(y)\\
& \le & \frac{\max_{y\in S^m\cap Z^\perp}|h(y)|}{\H^{m-1}(\S^{m-1})}
\int_{y\in S^m\cap Z^\perp}
F^{-m}(y)\,d\H^{m-1}(y)\\
& = & \left(\max_{y\in S^m\cap Z^\perp}|h(y)| \right)\cR[F^{-m}],
\end{eqnarray*}
so that we deduce from \eqref{AF-estimate} in connection with \eqref{AF-minkowski-bounds}
the gradient estimate for the Cartan area integrand at $Z\in S^m$
\begin{equation}\label{AF-estimate2}
|\cA^F_Z(Z)|^2 \le (\cA^F(Z))^2\Big(1+m^2 \sum_{\sigma=1}^{m+1}\left( \max_{y\in S^m\cap Z^\perp}|h(y)|^2\right)\Big)
\le M_F^{2m}\Big(1+(m+1)m^2\max_{y\in S^m\cap Z^\perp}|h(y)|^2\Big).
\end{equation}
To estimate $|h(y)|$ for $y\in S^m\cap Z^\perp$
recall that the fundamental tensor $(g_{ij})(y)=((F^2/2)_{y^iy^j})(y)$
is positive definite by (F2). Combining this with the homogeneity (F1) one can
easily deduce that the Hessian $F_{y^iy^j}(y)$ is positive semidefinite;
see, e.g., formula (1.2.9) in \cite[Proof of Theorem 1.2.2]{bao-chern-shen_2000}. Thus,
$$
\Lambda(F)|F_y(y)|^2\ge F_{y^i}(y)\Big(\frac{F^2}{2}\Big)_{y^iy^j}(y)F_{y^j}(y) = 
F_{y^i}(y)\big(F_{y^i}F_{y^j}+FF_{y^iy^j}\big)(y)F_{y^j}(y)
 \ge  |F_y(y)|^4
$$
for every $y\in \S^m$, where 
$$
\Lambda(F):=\max_{\eta,\zeta\in \S^m}\eta\cdot \big(F^2/2\big)_{yy}(\zeta)\eta
$$
is the largest possible eigenvalue of the fundamental tensor $g_{ij}$ of
$F$ when restricted to the sphere $\S^m$. (Recall the summation convention from
our introduction described under (F2).)
Since $F$ is bounded from below by $m_F>0$ on $\S^m$ (see \eqref{F-minkowski-bounds}),
we obtain
$$
|h(y)|^2\le \frac{1}{m_F^2}|Z|^2\Lambda(F)\quad\Foa y\in\S^m.
$$
Inserting this into \eqref{AF-estimate2} gives
\begin{equation}\label{gradient-bound}
|\cA^F_Z(Z)|\le M_F^m\sqrt{1+\Lambda(F)(m+1)\frac{m^2}{m_F^2}}\quad\Foa Z\in\S^m,
\end{equation}
which according to \eqref{winkl-isop} leads to the following isoperimetric
inequality
\begin{equation}\label{prelim-isop}
\area_\mathscr{M}^F(X)\le \frac{R}{m}\frac{M_F^m}{m_F}\sqrt{1+\Lambda(F)
(m+1)\frac{m^2}{m_F^2}}
\int_{\partial\mathscr{M}}dS.
\end{equation}
Recall from our remarks preceding Theorem \ref{thm:isop} that
$\int_{\partial\mathscr{M}}dS_F=\area^F_{\partial\mathscr{M}}(X|_{\partial\mathscr{M}})$
which specializes to 
$\int_{\partial\mathscr{M}}dS=\area^E_{\partial\mathscr{M}}
(X|_{\partial\mathscr{M}})$ if the Finsler metric $F$ happens to be
Euclidean, i.e.,  $F(x,y)=E(y)=|y|.$ Now, Lemma \ref{lem:area-comparison}
applied to $\Omega\equiv\Sigma:=\partial\mathscr{M}$ with dimension $k:=m-1$
leads to the desired isoperimetric inequality, since it implies
$\int_{\partial\mathscr{M}}dS
\le m_F^{-(m-1)}\int_{\partial\mathscr{M}}dS_F.$
\qed
\heikodetail{

\bigskip

Herr Overath hatte argumentiert:

{\tt \yy Da wir evtl. das Hausdorff-Ma\ss\  
$\H^{m-1}(X(\partial\mathscr{M}))$ in \eqref{prelim-isop} durch $\area^{|\cdot|}_{\partial\mathscr{M}}(X|_{\partial\mathscr{M}})$ ersetzen ist Lemma~\ref{lem:area-comparison} auf die $(m-1)$-dimensionale Immersion $X|_{\partial\mathscr{M}}:\partial\mathscr{M}\rightarrow \R^{m+1}$ anwendbar. Damit folgt $\int_{\partial \mathscr{M}} dS=\area^{|\cdot|}_{\partial\mathscr{M}}(X|_{\partial\mathscr{M}}) \le m_F^{-(m-1)} \area^{F}_{\partial\mathscr{M}}(X|_{\partial\mathscr{M}})=m_F^{-(m-1)}\int_{\partial \mathscr{M}} dS_F$.}

In conclusion we get
\begin{equation*}
\area_\mathscr{M}^F(X)\le \frac{R}{m}\left(\frac{M_F}{m_F}\right)^m\sqrt{1+\Lambda(F)(m+1)\frac{m^2}{m_F^2}}
\int_{\partial \mathscr{M}} dS_F.
\end{equation*}

\bigskip
\yy\yy\yy

\bigskip

}
For part (ii) we benefit from our analysis in the proof of part (i), and
use \cite[Theorem 3.2]{winklmann-isop} to obtain similarly to \eqref{winkl-isop}
the preliminary Finsler area estimate
\begin{equation}\label{winkl-isop2}
\area^F_\mathscr{M}(X)\le \|\cA^F_Z\|_{C^0(\S^2,\R^{3})}\sum_{i=1}^{k}\big[
\frac{\mathscr{L}(\Gamma_i)^2}{4\pi}+\frac 12  \mathscr{L}(\Gamma_i) \dist(a,\Gamma_i)\Big],
\end{equation}
where $\mathscr{L}(\Gamma)=\int |\dot{\Gamma}|$ denotes the Euclidean
length of  curve $\Gamma\subset\R^3$, and $\dist(a,\Gamma)$ is the
Eulcidean distance of a point $a\in\R^3$ to the curve $\Gamma\subset\R^3$.
To obtain a right-hand side that is completely expressed in terms
of Finslerian quantities we use, on the one hand, the estimate \eqref{gradient-bound}
on the gradient of the Cartan area integrand, and, on the other hand,
the following simple argument comparing Euclidean length and distance
to the Finslerian ones: by virtue of \eqref{F-minkowski-bounds} we have
\begin{equation}\label{length-compare}
\mathscr{L}(\Gamma)=\int |\dot{\Gamma}|\overset{\eqref{F-minkowski-bounds}}{\le}
\frac{1}{m_F}\int F(\dot{\Gamma})=\frac{1}{m_F}\mathscr{L}^F(\Gamma)
\end{equation}
for any curve $\Gamma\in\R^3.$ Since $\dist(a,\Gamma)$ and $\dist_F(a,\Gamma)$
between a
point $a\in\R^3$ and a curve $\Gamma\subset\R^3$ may be expressed by minimizing
the respective length functional over all curves connecting $a$ with some
 point on $\Gamma$ the estimate \eqref{length-compare}
immediately implies also
\begin{equation}\label{dist-compare}
\dist(a,\Gamma)\le\frac{1}{m_F}\dist_F(a,\Gamma).
\end{equation}
Inserting \eqref{gradient-bound}, \eqref{length-compare}, and \eqref{dist-compare}
into \eqref{winkl-isop2} we obtain the desired isoperimetric inequality
\eqref{isop2} as stated.

Notice for part (iii) that $F=F(x,y)$ is not a Minkowski metric, since it depends
explicitly on $x\in\R^3.$ The only isoperimetric inequalities for Cartan
functionals  depending on the position vector as well 
we are aware of are those of Clarenz and the second author in 
\cite{clarenz-vdM-isop}. 

By Corollary \ref{cor:ellipticity} we know that the Cartan area integrand
$\cA^F(x,Z)$ generated by the Finsler metric $F=F(x,y)$ is (parametric)
elliptic in the sense of Definition \ref{def:ellipticity}, in particular
we have 
\begin{equation}\label{specified-ellipticity}
\lambda |\pi_{Z^\perp}(\xi)|^2\le|Z|\xi\cdot\cA^F_{ZZ}(x,Z)\xi\le\Lambda |\pi_{Z^\perp}(\xi)|^2
\quad\Foa Z\in\R^3\setminus\{0\}, x\in \overline{B_1(0)},\,\xi\in\R^{3}
\end{equation}
for some constants $0<\lambda\le\Lambda<\infty$,
which is just \eqref{parametric-ellipticity} specified to dimension $m=2$
and for $R_0=1.$ (The upper bound in \eqref{specified-ellipticity} just
follows from smoothness of $\cA^F$ on $\R^{3}\times (\R^{3}\setminus\{0\})$.) 
Now, \eqref{specified-ellipticity} is exactly
condition (E) in \cite[p. 618]{clarenz-vdM-isop}, and since the positive
$1$-homogeneity holds for $\cA^F$ in any case, we can apply \cite[Theorem 2]{clarenz-vdM-isop} to find for the Euclidean area $\area_B(X)$ of $X(B)$
\begin{equation}\label{clarenz-vdM-isop}
\area_B(X)\le R\frac{2C(\cA^F)\frac{\Lambda}{\lambda}\big[\int_\Gamma\kappa\,ds-2\pi\big]
+\sqrt{\frac{\Lambda}{\lambda}}\mathscr{L}(\Gamma)}{2-Rh_F},
\end{equation}
if 
\begin{equation}\label{hF}
h_F:=C(\cA^F)(1+\|\cA^F_{xZ}\|_{\infty}^2)+\frac 1\lambda \|\cA^F_{xZ}\|_{\infty}<2,
\end{equation}
where $\|\cA^F_{xZ}\|_{\infty}:=
\|\sum_{i=1}^3\cA^F_{x_iZ_i}(\cdot,\cdot)\|_{C^0(\overline{B_1(0)}\times \S^2)}.$
\heikodetail{

\bigskip

\tt \yy Also $\|\cA^F_{xZ}\|_{\infty}:= \max_{i,j\in \{1,\ldots,m+1\}} \sup_{(x,Z)\in C^0(\overline{B_1(0)}\times \S^2,\R^{3\times 3})}|\cA^F_{xZ}(x,Z)|$? \yy

\bigskip

BERECHTIGTE FRAGE, ICH HABE ES NACHGESCHLAGEN UND OBEN KORRIGIERT\xx\xx \yy\yy Einverstanden \yy\yy

\bigskip

}
Here, 
$$
C(\cA^F):=C_x(\cA^F)+C_Z(\cA^F)\Big(1+\frac{\Lambda}{2\lambda^3}\Big)
$$ 
is a constant introduced in \cite[p. 628]{clarenz-vdM-isop}, where
$C_x(\cA^F)$ and $C_Z(\cA^F)$ are bounds on the $x$-derivative, and the $Z$-derivative
of the expression
\begin{equation}\label{lAF}
l^{\cA^F}(x,Z):=\frac{\cA^F_{ZZ}(x,Z)}{\sqrt{\det\cA^F_{ZZ}(x,Z)|_{Z^\perp}}}+
\frac{Z}{|Z|}\otimes
\frac{Z}{|Z|},
\end{equation}
(cf. formulae (35),(38) that lead to Proposition 2.4 (ii) in \cite{clarenz-vdM-isop}).
We need to investigate all quantities that enter the definition \eqref{hF} of $h_F$,
since we need to verify its smallness.

Since the denominator in \eqref{lAF} is bounded by $\lambda$ from below by virtue of
\eqref{specified-ellipticity} we find that for each fixed $x\in\overline{B_1(0)}\subset
\R^3$
$$
\varrho_1(l^{\cA^F}(x,\cdot)-l^{\cA^E}(\cdot))\longrightarrow 0\quad\textnormal{as \,\,$
\varrho_3(\cA^F(x,\cdot)-\cA^E(\cdot))\to 0$,}
$$
where we refer to the definition of the spherical seminorms in \eqref{semi-norms-sphere}
and to our notation of the Euclidean metric $E(y)=|y|$ on $\R^3.$
Analogously to the proof of part (ii) of Theorem \ref{thm:bernstein} we can show
that for given $\epsilon >0$ we find $\delta=\delta(\epsilon)$ such that
the inequality \eqref{C3-Finsler-close} implies
that
\begin{equation}\label{eps-criterion}
\varrho_1(l^{\cA^F}(x,\cdot)-l^{\cA^E}(\cdot))<\epsilon\quad\Foa x\in\overline{B_1(0)}.
\end{equation}
One easily checks in \eqref{lAF} that $l^{\cA^E}(Z)=\Id_{\R^3}$ so that
$l^{\cA^E}_Z=0$, and we obviously have no $x$-dependence in $l^{\cA^E}$, hence
$l^{\cA^E}_x=0$. In addition, since $\cA^E$ does not explicitly depend on $x$ we
have $\cA^E_{xZ}=0$. 
\heikodetail{

With 
$$
\cA^E_{ZZ}=\frac{\Id}{|Z|}-\frac{Z\otimes Z}{|Z|^3}
$$
one finds
$$
\cA^E_{ZZ}Z=0\quad\AND\quad \cA^E_{ZZ}(Z)\eta=\frac{1}{|Z|}\eta\Foa\eta\in Z^\perp,
$$
which implies that $\lambda_1=0$ is a simple eigenvalue of $cA^E_{ZZ}$ and
$\lambda_2=1/|Z|$ is an eigenvalue of $cA^E_{ZZ}$ with algebraic multiplicity
$2$.
Hence $\det cA^E_{ZZ}(Z)|_{Z^\perp}=1/|Z|^2$, which implies that
$$
l^{\cA^E}=|Z|\Big[\frac{\Id}{|Z|}-\frac{Z\otimes Z}{|Z|^3}\Big]+\frac{Z}{|Z|}\otimes
\frac{Z}{|Z|}=\Id.
$$

\bigskip

}
Consequently, by \eqref{eps-criterion}, we can choose $\delta_1$
sufficiently small such that for all $\delta\in (0,\delta_1)$ for which 
\eqref{C3-Finsler-close} holds , we have
$$
\|\cA^F_{xZ}\|^2_\infty<\min\{1,\lambda\}\quad\AND\quad C(\cA^F)(1+
\|\cA^F_{xZ}\|^2_\infty)<1,
$$
so that \eqref{hF} is satisfied. 
\heikodetail{

\bigskip

\yy Notice that $\cA^E_{xZ}\equiv 0$.\yy

\bigskip

PLEASE EXPLAIN\xx\xx \yy\yy Die Frage hat sich erledigt, die Identitaet wurde kurz vorher im Beweis erwaehnt! \yy\yy

\bigskip

} 
Therefore, we have indeed the preliminary
isoperimetric inequality \eqref{clarenz-vdM-isop}. Since we assumed
 $R\le 1$,
\heikodetail{

\bigskip

\tt\yy This is only the case if we set $R$ to the special choice given in \cite{clarenz-vdM-isop}, namely $R:=\inf_{q\in\R^3}\|X(\cdot)-q \|_{\infty,\partial B}$? \yy

\bigskip

BERECHTIGTER EINWAND, HABE DIE VORAUSSETZUNGEN VERSTAERKT UND OBEN KORRIGIERT,
die Voraussetzung $\|X\|_{L^\infty,\R^3}\le 1$ ist hier direkt nicht mehr sichtbar eingegangen, wird aber in \eqref{clarenz-vdM-isop} gebraucht, und zwar m.E.
nicht nur, um das dortige $R_\Gamma=\inf_{q\in\R^3}\|X(\cdot)-q \|_{\infty,\partial B}$ nach oben gegen $1$ abzuschaetzen, oder sind Sie da anderer Meinung?\xx\xx \yy\yy Ich bin jetzt ganz Ihrer Meinung!\yy\yy

\bigskip

}
we can omit $R$ in the denominator. The simple
estimates
$$
\area^F_B(X)\le {M^*_F}^2\area_B(X)\quad\AND\quad \mathscr{L}(\Gamma)\le\mathscr{L}^F(\Gamma)/m_F^*
$$
(cf.  Lemma \ref{lem:area-comparison} and \eqref{length-compare}) lead to the desired isoperimetric
inequality \eqref{isop3}, 
where we have set 
$$
c_1(F):=\frac{2C(\cA^F)\frac{\Lambda}{\lambda}}{2-h_F}\quad\AND\quad
c_2(F):=\frac{\sqrt{\Lambda/\lambda}}{2-h_F},
$$
which also explains why $c_1$ vanishes  and $c_2=1/2$
 if $F=E(y)=|y|$  since $C(\cA^E)=0$ and $h_F=h_E=0$ verifying our remarks
 following Theorem \ref{thm:isop}.
 \qed

\bigskip

\noindent
{\bf Existence and uniqueness of Finsler-minimal immersions spanning given
boundary contours.}\,
Up to now there are only very few results guaranteeing the existence of
immersed surfaces minimizing a general Cartan functional and spanning
a given boundary contour. We will draw from White's existence result
under the condition of  
extreme boundary curves \cite{white-existence-1991} to prove 
with our results of Section \ref{sec:2} the 
existence of embedded Finsler-minimal disks spanning such a given
boundary curve in general Finsler spaces
with metric $F=F(x,y)$. In the Minkowski case, $F=F(y)$, we will
benefit from  a simpler  existence 
and uniqueness result of Winklmann \cite{winklmann-existence-uniqueness}.

Let us first explain the notion of {\it Finsler area mean convexity}
in our assumption in part (i) of Theorem \ref{thm:existence} 
adopting the variational characterization that White \cite{white-existence-1991} 
used for 
Cartan functionals but without using Finsler-mean curvature
in our context to avoid the ambivalent choice of a suitable Finsler
normal; see also \cite{bergner-froehlich-existence-2009} for an
equivalent notion for
Cartan integrands that do not depend on position.
\begin{definition}[\textbf{Inward-variations and Finsler area mean convexity}]\label{def:inward-mean-convex}
Let $Y:\partial\Sigma\to\R^3$ be the injection of the boundary $\partial\Sigma$
of a smooth $3$-dimensional submanifold $\Sigma\subset\R^3$ 
into $\R^3$. 
\begin{enumerate}
\item[\rm (i)]
A
 smooth mapping $\tilde{Y}:(-\eps,\eps)\times\partial\Sigma\to\R^3$ for some
 $\eps>0$
is called an {\em inward variation of $\partial\Sigma$} if $\tilde{Y}(0,\cdot)=Y(\cdot)$ on $\partial\Sigma$, with
a smooth variation vector field 
 $V(\cdot):=\frac{d}{dt}|_{t=0}\tilde{Y}(t,\cdot):\partial\Sigma\to\R^3$, not identically zero, that satisfies for every $u\in\partial\Sigma$   either $V(u)=0$ or 
$V(u)\not\in T_u\partial\Sigma$ and $V(u)=c'(s)$ for some smooth curve
$c:[0,\delta)\to\Sigma$, $\delta>0$, with $c(0)=u$.
\item[\rm (ii)]
If $\R^3$ is equipped with a Finsler metric $F=F(x,y)$ then $\Sigma\subset\R^3$
is said to be {\em strictly Finsler area mean convex} if
$$
\frac{d}{dt}|_{t=0}\area^F_{\partial\Sigma}(\tilde{Y}(t,\cdot))<0
$$
for all inward variations of $\partial\Sigma$.
\end{enumerate}
\end{definition}
{\sc Proof of Theorem \ref{thm:existence}.}\,
The explicit representation of Finsler area  \eqref{area_space} in
Theorem \ref{thm:finsler-area}  implies that a strictly Finsler
area mean convex $3$-dimensional submanifold $\Sigma\subset\R^3$
with boundary $\partial\Sigma$ is strictly $\cA^F$-convex in the 
sense of White \cite[items 1.3 \& 1.5]{white-existence-1991}. (Notice
that our smoothness assumption on $\Sigma $ implies the strict pointwise
inequality required for the normal component of  White's $\cA^F$-mean curvature
 \cite[1.5]{white-existence-1991}
even everywhere, not only almost everywhere. White's more general class
of Lipschitz variations (see \cite[1.2]{white-existence-1991}), on the
other hand, lead to the same variational equations for Cartan functionals
as, e.g.,  in the work of Clarenz et al. \cite{clarenz-phd_1999,clarenz_enclosure,clarenz-vdM-compact_2001}, and  
therefore to Finsler minimal surfaces in our context.)
White's ellipticity assumption \cite[pp.  413,414]{white-existence-1991} on the  
Cartan integrand $\cA^F$ reads as the uniform
convexity of the $\cA^F(x,\cdot)$-unit ball $\{Z\in\R^3:\cA^F(x,Z)\le 1\}$
for each $x$ in $\R^3,$ which according to \cite[p. 72]{gia-hil-II_1996}
is 
satisfied if $\cA^F $ is parametric elliptic in the sense of 
Definition \ref{def:ellipticity}.  That this is indeed the case is stated
in Corollary \ref{cor:ellipticity}. 
 Hence White's existence
result \cite[version of Theorem 3.4 on p. 425]{white-existence-1991} is applicable
in our context, which in particular leads to the first statement
in part (i) of Theorem \ref{thm:existence}. 
\heikodetail{

\bigskip
\tt
ACHTUNG: Sollen wir noch etwas zu dem Fall vorgeschriebenen Geschlechts
sagen, siehe White's allgemeine Version von Thm 3.4 auf Seite 13!\xx\xx \yy\yy In der Doktorarbeit habe ich White's genus-Formulierung mit aufgenommen. Hier koennte man das auch, ich denke aber, dass es auch ohne funktioniert. \yy\yy

\bigskip

}

The second statement in part (i) is stated in \cite[version of Theorem 3.4 on p. 413]{white-existence-1991} and 
follows from the fact that for a Minkowski metric $F=F(y)$ on $\R^3$
smooth convex bodies $\Sigma\subset\R^3$ are Finsler area mean convex; see
also the other examples of Bergner and Fr\"ohlich in \cite[p. 368]{bergner-froehlich-existence-2009} whose concept of weighted mean convexity coincides with the 
one of White 
in case of Cartan integrands without $x$-dependence.

For part (ii) we refer to \cite[Corollary 1.3]{winklmann-existence-uniqueness},
where we have to restrict to a Minkowski metric $F=F(y)$ since then the
corresponding Cartan area functional $\cA^F$ depends on $Z$ only, so that
the Finsler area 
$$
\area^F_\Omega(X)=\int_\Omega\cA^F(X_{u^1}\wedge X_{u^2})\,du^1\wedge du^2, \quad\Omega\subset\R^2,
$$ 
belongs to the Cartan functionals considered in \cite{winklmann-existence-uniqueness}.
The relevant homogeneity and ellipticity conditions formulated in
\cite[p. 269]{winklmann-existence-uniqueness} are satisfied by $\cA^F$ as we have
observed before; see, e.g., \eqref{AF-minkowski-ellipticity}, and
Finsler-minimal immersions are exactly  {\it $\cA^F$-minimal immersions} in
the language of  Clarenz et al. and Winklmann, which means that it has vanishing $\cA^F$-mean
curvature $H_{\cA^F}=0$; see \cite[p. 270]{winklmann-existence-uniqueness}.

\addcontentsline{toc}{section}{References}

\bibliography{finsler-plateau}{}

\begin{thebibliography}{10}

\bibitem{alvarez-berck-wrong-hausdorff-2006}
J.~C. {\'A}lvarez~Paiva and G.~Berck.
\newblock What is wrong with the {H}ausdorff measure in {F}insler spaces.
\newblock {\em Adv. Math.}, 204(2):647--663, 2006.

\bibitem{bailey}
T.~N. Bailey, M.~G. Eastwood, A.~R. Gover, and L.~J. Mason.
\newblock Complex analysis and the {F}unk transform.
\newblock {\em J. Korean Math. Soc.}, 40(4):577--593, 2003.
\newblock Sixth International Conference on Several Complex Variables
  (Gyeongju, 2002).

\bibitem{bao-chern-shen_2000}
D.~Bao, S.-S. Chern, and Z.~Shen.
\newblock {\em An introduction to {R}iemann-{F}insler geometry}, volume 200 of
  {\em Graduate Texts in Mathematics}.
\newblock Springer-Verlag, New York, 2000.

\bibitem{bergner-froehlich-existence-2009}
Matthias Bergner and Steffen Fr{\"o}hlich.
\newblock Existence, uniqueness and graph representation of weighted minimal
  hypersurfaces.
\newblock {\em Ann. Global Anal. Geom.}, 36(4):363--373, 2009.

\bibitem{busemann-1947}
Herbert Busemann.
\newblock Intrinsic area.
\newblock {\em Ann. of Math. (2)}, 48:234--267, 1947.

\bibitem{busemann-convex-brunn-minkowski-1949}
Herbert Busemann.
\newblock A theorem on convex bodies of the {B}runn-{M}inkowski type.
\newblock {\em Proc. Nat. Acad. Sci. U. S. A.}, 35:27--31, 1949.

\bibitem{chern-shen-2005}
Shiing-Shen Chern and Zhongmin Shen.
\newblock {\em Riemann-{F}insler geometry}, volume~6 of {\em Nankai Tracts in
  Mathematics}.
\newblock World Scientific Publishing Co. Pte. Ltd., Hackensack, NJ, 2005.

\bibitem{clarenz-phd_1999}
Ulrich Clarenz.
\newblock {\em S\"atze \"uber Extremalen zu parametrischen Funktionalen}.
\newblock PhD thesis, Univ. of Bonn, 1999.

\bibitem{clarenz_enclosure}
Ulrich Clarenz.
\newblock Enclosure theorems for extremals of elliptic parametric functionals.
\newblock {\em Calc. Var. Partial Differential Equations}, 15(3):313--324,
  2002.

\bibitem{clarenz-vdM-compact_2001}
Ulrich Clarenz and Heiko von~der Mosel.
\newblock Compactness theorems and an isoperimetric inequality for critical
  points of elliptic parametric functionals.
\newblock {\em Calc. Var. Partial Differential Equations}, 12(1):85--107, 2001.

\bibitem{clarenz-vdM-isop}
Ulrich Clarenz and Heiko von~der Mosel.
\newblock Isoperimetric inequalities for parametric variational problems.
\newblock {\em Ann. Inst. H. Poincar\'e Anal. Non Lin\'eaire}, 19(5):617--629,
  2002.

\bibitem{cui-shen-2009}
Ningwei Cui and Yi-Bing Shen.
\newblock Bernstein type theorems for minimal surfaces in
  {$(\alpha,\beta)$}-space.
\newblock {\em Publ. Math. Debrecen}, 74(3-4):383--400, 2009.

\bibitem{cui-shen-2011}
Ningwei Cui and Yi-Bing Shen.
\newblock Minimal rotational hypersurfaces in {M}inkowski
  {$(\alpha,\beta)$}-space.
\newblock {\em Geom. Dedicata}, 151:27--39, 2011.

\bibitem{DHKW1}
Ulrich Dierkes, Stefan Hildebrandt, Albrecht K{\"u}ster, and Ortwin Wohlrab.
\newblock {\em Minimal surfaces. {I}}, volume 295 of {\em Grundlehren der
  Mathematischen Wissenschaften [Fundamental Principles of Mathematical
  Sciences]}.
\newblock Springer-Verlag, Berlin, 1992.
\newblock Boundary value problems.

\bibitem{evans-gariepy}
Lawrence~C. Evans and Ronald~F. Gariepy.
\newblock {\em Measure theory and fine properties of functions}.
\newblock Studies in Advanced Mathematics. CRC Press, Boca Raton, FL, 1992.

\bibitem{funk-1915}
Paul Funk.
\newblock \"{U}ber eine geometrische {A}nwendung der {A}belschen
  {I}ntegralgleichung.
\newblock {\em Math. Ann.}, 77(1):129--135, 1915.

\bibitem{gardner}
Richard~J. Gardner.
\newblock {\em Geometric tomography}, volume~58 of {\em Encyclopedia of
  Mathematics and its Applications}.
\newblock Cambridge University Press, Cambridge, 1995.

\bibitem{gia-hil-II_1996}
Mariano Giaquinta and Stefan Hildebrandt.
\newblock {\em Calculus of variations. {II}}, volume 311 of {\em Grundlehren
  der Mathematischen Wissenschaften [Fundamental Principles of Mathematical
  Sciences]}.
\newblock Springer-Verlag, Berlin, 1996.
\newblock The Hamiltonian formalism.

\bibitem{gilbarg-trudinger_1998}
David Gilbarg and Neil~S. Trudinger.
\newblock {\em Elliptic partial differential equations of second order}.
\newblock Classics in Mathematics. Springer-Verlag, Berlin, 2001.
\newblock Reprint of the 1998 edition.

\bibitem{groemer}
H.~Groemer.
\newblock {\em Geometric applications of {F}ourier series and spherical
  harmonics}, volume~61 of {\em Encyclopedia of Mathematics and its
  Applications}.
\newblock Cambridge University Press, Cambridge, 1996.

\bibitem{helgason_radon}
Sigurdur Helgason.
\newblock {\em The {R}adon transform}, volume~5 of {\em Progress in
  Mathematics}.
\newblock Birkh\"auser Boston, Mass., 1980.

\bibitem{helgason_groups}
Sigurdur Helgason.
\newblock {\em Groups and geometric analysis}, volume 113 of {\em Pure and
  Applied Mathematics}.
\newblock Academic Press Inc., Orlando, FL, 1984.
\newblock Integral geometry, invariant differential operators, and spherical
  functions.

\bibitem{Hil-sauvigny-energy-estimate-2009}
Stefan Hildebrandt and Friedrich Sauvigny.
\newblock An energy estimate for the difference of solutions for the
  {$n$}-dimensional equation with prescribed mean curvature and removable
  singularities.
\newblock {\em Analysis (Munich)}, 29(2):141--154, 2009.

\bibitem{HilvdM-dominance}
Stefan Hildebrandt and Heiko von~der Mosel.
\newblock Dominance functions for parametric {L}agrangians.
\newblock In {\em Geometric analysis and nonlinear partial differential
  equations}, pages 297--326. Springer, Berlin, 2003.

\bibitem{HilvdM-parma}
Stefan Hildebrandt and Heiko von~der Mosel.
\newblock Conformal representation of surfaces, and {P}lateau's problem for
  {C}artan functionals.
\newblock {\em Riv. Mat. Univ. Parma (7)}, 4*:1--43, 2005.

\bibitem{jenkins_1961b}
H.~B. Jenkins.
\newblock On two-dimensional variational problems in parametric form.
\newblock {\em Arch. Rational. Mech. Anal.}, 8:181--206, 1961.

\bibitem{overath-phd_2013}
Patrick Overath.
\newblock {\em Minimal immersions in Finsler spaces}.
\newblock PhD thesis, RWTH Aachen University, Feb. 2014.

\bibitem{overath-vdM-2012a}
Patrick Overath and Heiko von~der Mosel.
\newblock Plateau's problem in {F}insler $3$-space.
\newblock {\em Manuscripta Mathematica}, 143:273--316, 2014.

\bibitem{radon-1994}
Johann Radon.
\newblock \"{U}ber die {B}estimmung von {F}unktionen durch ihre {I}ntegralwerte
  l\"angs gewisser {M}annigfaltigkeiten.
\newblock In {\em 75 years of {R}adon transform ({V}ienna, 1992)}, Conf. Proc.
  Lecture Notes Math. Phys., IV, pages 324--339. Int. Press, Cambridge, MA,
  1994.

\bibitem{raewer-phd_1993}
Klaus R{\"a}wer.
\newblock {\em Stabile Extremalen parametrischer Doppelintegrale in {$\R^3$}}.
\newblock PhD thesis, Univ. of Bonn, 1993.

\bibitem{rudin-fa-book_1973}
Walter Rudin.
\newblock {\em Functional analysis}.
\newblock McGraw-Hill Book Co., New York, 1973.
\newblock McGraw-Hill Series in Higher Mathematics.

\bibitem{shen98}
Zhongmin Shen.
\newblock On {F}insler geometry of submanifolds.
\newblock {\em Math. Ann.}, 311(3):549--576, 1998.

\bibitem{simon_1977b}
Leon Simon.
\newblock On a theorem of de {G}iorgi and {S}tampacchia.
\newblock {\em Math. Z.}, 155(2):199--204, 1977.

\bibitem{simon_1977a}
Leon Simon.
\newblock On some extensions of {B}ernstein's theorem.
\newblock {\em Math. Z.}, 154(3):265--273, 1977.

\bibitem{sst}
Marcelo Souza, Joel Spruck, and Keti Tenenblat.
\newblock A {B}ernstein type theorem on a {R}anders space.
\newblock {\em Math. Ann.}, 329(2):291--305, 2004.

\bibitem{souza-tenenblat-2003}
Marcelo Souza and Keti Tenenblat.
\newblock Minimal surfaces of rotation in {F}insler space with a {R}anders
  metric.
\newblock {\em Math. Ann.}, 325(4):625--642, 2003.

\bibitem{white-existence-1991}
Brian White.
\newblock Existence of smooth embedded surfaces of prescribed genus that
  minimize parametric even elliptic functionals on {$3$}-manifolds.
\newblock {\em J. Differential Geom.}, 33(2):413--443, 1991.

\bibitem{winklmann-isop}
Sven Winklmann.
\newblock Isoperimetric inequalites involving generalized mean curvature.
\newblock {\em Analysis (Munich)}, 22(4):393--403, 2002.

\bibitem{winklmann-existence-uniqueness}
Sven Winklmann.
\newblock Existence and uniqueness of {$F$}-minimal surfaces.
\newblock {\em Ann. Global Anal. Geom.}, 24(3):269--277, 2003.

\bibitem{winklmann-bernstein}
Sven Winklmann.
\newblock A {B}ernstein result for entire {$F$}-minimal graphs.
\newblock {\em Analysis (Munich)}, 27(4):375--386, 2007.

\end{thebibliography}
\bibliographystyle{plain}


\small
\vspace{1cm}
\begin{minipage}{56mm}
{\sc Patrick Overath}\\
Institut f\"ur Mathematik\\
RWTH Aachen University\\
Templergraben 55\\
D-52062 Aachen\\
GERMANY\\
E-mail: {\tt overath@}\\
{\tt instmath.rwth-aachen.de}
\end{minipage}
\hfill
\begin{minipage}{56mm}
{\sc Heiko von der Mosel}\\
Institut f\"ur Mathematik\\
RWTH Aachen University\\
Templergraben 55\\
D-52062 Aachen\\
GERMANY\\
Email: {\tt heiko@}\\{\tt instmath.rwth-aachen.de}
\end{minipage}

\end{document}